\documentclass[11pt]{article}

\usepackage[pdftex]{graphicx} 
\usepackage{textcomp}
\usepackage{amsbsy}
\usepackage{latexsym}
\usepackage[mathscr]{eucal}
\usepackage{amsfonts,amsmath,amsthm}
\usepackage{amssymb}
\usepackage[dvipsnames]{xcolor}
\usepackage{graphicx}
\usepackage{enumitem}
\usepackage{array, multirow, multicol}

\usepackage{algorithm}
\usepackage[noend]{algpseudocode}

\makeatletter
\def\BState{\State\hskip-\ALG@thistlm}
\makeatother

\providecommand{\classification}[1]{\noindent\textbf{Mathematics Subject Classification~} #1}

\usepackage{fullpage}

\definecolor{aggie}{cmyk}{0.15,1.0,0.39,0.69}

\newtheorem{lemma}{Lemma}[section]

\newtheorem{cor}[lemma]{Corollary}
\newtheorem{theorem}[lemma]{Theorem}
\newtheorem{remark}[lemma]{Remark}

\def\RR{\rm \hbox{I\kern-.2em\hbox{R}}}
\def\NN{\rm \hbox{I\kern-.2em\hbox{N}}}
\def\ZZ{\rm {{\rm Z}\kern-.28em{\rm Z}}}
\def\CC{\rm \hbox{C\kern -.5em {\raise .32ex \hbox{$\scriptscriptstyle
|$}}\kern
-.22em{\raise .6ex \hbox{$\scriptscriptstyle |$}}\kern .4em}}

\def\<{\langle}
\def\>{\rangle}

\def\e{\varepsilon}

\def\cB{{\cal B}}
\def\cF{{\cal F}}
\def\cU{{\cal U}}

\def\cP{{\cal P}}

\def\cO{{\cal O}}
\def\Chi{\raise .3ex
\hbox{\large $\chi$}} 
\def\lsima{\hbox{\kern -.6em\raisebox{-1ex}{$~\stackrel{\textstyle<}{\sim}~$}}\kern -.4em}
\def\lsim{\hbox{\kern -.2em\raisebox{-1ex}{$~\stackrel{\textstyle<}{\sim}~$}}\kern -.2em}

\def\R{\mathbb{R}}

\def\T{{\relax\ifmmode I\!\!\hspace{-1pt}T\else$I\!\!\hspace{-1pt}T$\fi}}
\def\N{\mathbb{N}}

\def\N{\mathbb{N}}

\def\C{\mathbb{C}}

\def\lsim{\raisebox{-1ex}{$~\stackrel{\textstyle<}{\sim}~$}}

  \def\NN{N}                  

\def\cU{{\cal U}}

\def\cO{{\cal O}}

\def\cU{{\cal U}}

\def\cL{{\cal L}}
\def\cB{{\cal B}}

\def\cP{{\cal P}}
\def\cQ{{\cal Q}}

\def\cU{{\cal U}}

\def \meas {{\rm meas}}

\newcommand{\be}{\begin{equation}}
\newcommand{\ee}{\end{equation}}
\newcommand{\bes}{\begin{equation*}}
\newcommand{\ees}{\end{equation*}}
\newcommand{\bea}{$$ \begin{array}{lll}}
\newcommand{\eea}{\end{array} $$}

\def \exp{\mathop{\rm    exp}}

\newcommand{\beqn}{\begin{equation}}
\newcommand{\eeqn}{\end{equation}}

\def\endproof{\hfill\rule{1.5mm}{1.5mm}\\[2mm]}

\makeatletter
\@addtoreset{equation}{section}
\makeatother

\newcommand\eref[1]{{\rm (\ref{#1})}}
\newcommand{\iref}[1]{{\rm (\ref{#1})}}

\def\int{\intop\limits}

\newcommand\cH{{\cal H}}

\title
{ Polynomial Approximation of  Anisotropic Analytic Functions of Several Variables}
\author{ 
 Andrea Bonito, Ronald DeVore, Diane Guignard, Peter Jantsch, and Guergana Petrova
\thanks{
This research was supported by the NSF grants DMS-1817691 (AB),   DMS 15-21067 (RD-GP), DMS 18-17603 (RD-GP), ONR grants N00014-17-1-2908 (RD), N00014-16-1-2706 (RD); DG was supported by the Swiss National Science Foundation grant P2ELP2-175056 and  IAMCS at TAMU, and PJ was supported by an NSF Fellowship DMS-1704121. A portion of this research was completed while RD (Simon Fellow), DG, and PJ were supported as  visitors of the Isaac Newton Institute at Cambridge University.}}

\begin{document}

\maketitle
\date{}

\begin{abstract}   
Motivated by numerical methods for solving parametric partial differential equations, this paper studies the approximation of multivariate analytic functions by algebraic polynomials. We introduce various anisotropic  model  classes based on Taylor expansions, and study their approximation by finite dimensional polynomial spaces $\cP_\Lambda$ described by lower sets $\Lambda$.  Given a budget $n$ for the dimension of $\cP_\Lambda$, we prove that certain lower sets  $\Lambda_n$, with cardinality $n$, provide a  certifiable approximation error that is in a certain sense optimal, and that these lower sets have a simple definition in terms of simplices. Our main goal is to obtain approximation results  when  the number of variables $d$ is large and even infinite, and so we concentrate almost exclusively on the case $d=\infty$. We also emphasize obtaining results which hold for the full range $n\ge 1$, rather than asymptotic results that only hold for $n$ sufficiently large. In applications, one typically wants $n$ small to comply with computational budgets.
\end{abstract}
\smallskip

\classification{ 41A10, 41A58, 41A63, 65N15}

\section{Introduction} 
 
Polynomial and piecewise polynomial approximation are a staple in numerical analysis. For example, approximation by piecewise polynomials on simplicial partitions is the underpinning of Finite Element Methods. In that setting, one approximates the solution $u$ to a partial differential equation (PDE) on a domain $D\subset \R^{ d}$, where $d$ is typically small ($d=1,2,3$). The solution $u$ to the PDE typically has limited regularity, and the rate of approximation is of order $O(n^{-r})$, where $n$ is the number of degrees of freedom  in the approximation and $r$ is small. This type of approximation is well-understood by means of theorems which  relate the approximation order  $r$ to the smoothness order $s$ of $u$ in certain Sobolev and Besov spaces (see \cite{BDD,BDDP, BDN, BS}). The approximation rate takes the form $r=s/d$ and therefore deteriorates as  $d$ increases. This is commonly referred to as the curse of dimensionality.

The present paper is interested in a different setting that  arises in other application areas, in particular when using  numerical methods for solving  stochastic or {\em parametric} PDEs. In that setting, one wishes to approximate the solution $u$  to the parametric PDE which depends on input parameters $y$ and takes values in a Banach space $X$. The parameters $y$ come from a set $Y\subset \R^d$ where $d$ is large or even infinite. Hence, it is often crucial to perform a model reduction (dimension reduction) for the solution map $y \rightarrow u(y) \in X$ of the  parametric PDE. One possibility to obtain such dimension reduction is to approximate $u$ by Banach space valued polynomials in $y$. The main property of $u$ that makes such an approximation possible is that under standard  assumptions on the parametrized coefficients of the PDE, it is known that $u$ admits an analytic extension onto certain complex polydiscs that contain $Y$ (see \cite{CD}). In other words, $u$ has a certain anisotropic analyticity. This motivates the study of approximation of anisotropic analytic functions by polynomials, which is the subject of the present  paper. Although we are motivated by parametric PDE applications, we formulate and study this subject as purely a problem in multivariate approximation. In this way, we hope to draw the attention of the approximation community to this area of research.

For the most part, we are interested in the case  of an infinite number of parameters, i.e., $d=\infty$.  This allows us to prove results which are immune to the dimension $d$ and is a common setting in parametric PDEs. Specifically, we take parameters in the set $Y:=[-1,1]^\N$, where $\N$ is the set of natural numbers. Sometimes we remark on the   case $ Y_d:=[-1,1]^d$ with $d$ finite, in particular, when we wish to compare our results with other results in the literature established only for finite $d$.   

Let $\cF$ denote the set of all infinite sequences $\nu=(\nu_1,\nu_2,\dots)$ with entries  $\nu_j\in\N_0:=\N\cup \{0\}$, where only a finite number of the entries in $\nu$  are allowed to be nonzero. If $\Lambda\subset \cF$ is a finite subset of $\cF$, we denote by $\cP_\Lambda$  the space of $X$-valued polynomials spanned by the monomials $y^\nu$, where the $\nu$ come from the set $\Lambda$. Thus, any element of $\cP_\Lambda$ has the form

\be
\label{polynomials} 
P(y)=\sum_{\nu\in \Lambda} c_\nu y^\nu,
\ee
where the coefficients $c_\nu$ come from $X$. Here and throughout the paper, we use standard multivariate notation.  In particular, $y^\nu:=y_1^{\nu_1}y_2^{\nu_2}\cdots$. Since $\nu$ has only a finite number of nonzero entries, any such product is finite. 

For any $u\in L_\infty(Y,X) $, and any finite set $\Lambda\subset \cF$, we define the error of approximation of $u$ by polynomials in $\cP_\Lambda$ to be
\be
\label{aerror}
E_\Lambda(u):=\inf_{P\in\cP_\Lambda}\|u-P\|_{L_\infty(Y,X)}, \quad \hbox{and}\quad \|v\|_{L_\infty(Y,X)}:=\sup_{y\in Y}\|v(y)\|_X,
\ee
where $L_\infty(Y,X)$ consists of all functions $v$ on $Y$ that are bounded mappings into $X$.

If no conditions are imposed, then the potential sets $\Lambda$ may be quite complex and beyond the scope of numerical methods. For this reason, one usually imposes additional structure on these sets such as fixed total degree or fixed coordinate degree in the case $d$ is finite. We are especially interested in the case where the sets $\Lambda$ are {\it lower sets}, that is, sets $\Lambda$ with the property
$$
\hbox {if} \quad \nu\in \Lambda, \quad \hbox{then} \quad \mu\in\Lambda \quad \hbox{whenever} \quad \mu_j\le \nu_j, \quad j=1,2, \dots.
$$
We consider the collection $\cL_n$ of lower sets with cardinality $\leq n$,
$$
\cL_0:=\emptyset, \quad \cL_n:=\{\Lambda \subset \cF:\,\,\#\Lambda \leq n,\,\,\Lambda \  \hbox{is a lower set}\}, \quad n=1,2,\dots.
$$
and for given compact class $K$ of functions in $L_\infty(Y,X)$, and any finite set $\Lambda\subset \cF$, we define
\be
\label{classerror}
E_\Lambda(K):= \sup_{u\in K}E_\Lambda(u),
\ee
and 
\be
\label{classerrorlower}
E_0(K):=\sup_{u\in K}\|u\|_{L_\infty(Y,X)}, \quad E_n(K):=\inf_{\Lambda \in\cL_n}E_\Lambda(K), \quad n=1, 2,\dots.
\ee
Notice that in this definition the set $\Lambda$ is allowed to depend on $K$, but cannot change for the various $u\in K$.  So, as formulated,  this is a problem
of finding the best linear space $\cP_\Lambda$ to use when approximating $K$. Hence, the optimal performance $E_n(K)$ satisfies
\be
\label{comparewidths}
d_n(K)_{L_\infty(Y,X)}\le E_n(K),
\ee
where $d_n$ denotes the Kolmogorov $n$-width of $K$ in $L_\infty(Y,X)$. The sets $K$ are commonly called {\it model classes}. The case where the error of approximation is measured in $L_q(Y,X)$, $q<\infty$, is also interesting but not studied here.

Given the  model class $K$, we are interested in several fundamental issues: 
\begin{itemize}
\item First, what can we say about the  rate of decay of $E_n(K)$ as $n$ increases? By now, there are several results in the literature that give upper bounds on  $E_n(K)$ for certain anisotropic analytic classes $K$ of the type analyzed in this paper. Most often these bounds have been developed in the setting where $u$ is a solution to a parametric PDE. So part of our effort is to separate out which of these results are simply a result of the analyticity of $u$ and do not use any additional  properties of the PDE solution.

\item A second important issue is the optimality of the known bounds for $E_n(K)$. Indeed, the typical results only give upper bounds for $E_n(K)$ and sometimes only for $n$ sufficiently large.

\item  A third significant problem  in this area of research is to give a recipe for finding good lower sets $\Lambda_n\in\cL_n$ such that  $E_{\Lambda_n}(K)$ performs at or near   $E_n(K)$.  This can be  a nontrivial issue in numerical applications since, given a budget $n$, searching over all lower sets in $\cL_n$   to find a suitable $\Lambda_n$ is prohibitive.

\end{itemize}

As already noted, we are interested in model classes $K$ described by some form of anisotropic analyticity.   We focus on $X$ valued functions that are  analytic on a polydisc $D_\rho$ consisting of all complex sequences $z=(z_1, z_2,\dots)$, with $|z_j|<\rho_j$, $j=1,2,\dots$. Here, $\rho=(\rho_1,\rho_2,\dots)$ is always a nondecreasing sequence of positive numbers with $\rho_1>1$. The functions in $D_\rho$ have more smoothness in the variable $z_j$ as $j$ increases. In turn, the influence of this variable on the value of $u$ at a point in $Y$ is weaker. Any function in $D_\rho$ has Taylor coefficients that are elements of $X$. In \S \ref{sec:anisoclass}, we introduce a variety of spaces $\cB_{\rho,p}$, which differ in the assumptions imposed on the Taylor coefficients. These spaces are motivated by recent work (see \cite{BCM,CM,CD}) in parametric PDEs.
 
The remainder of the paper concentrates on understanding the rate of decay of $E_n(K)$ for these model classes and understanding how to choose lower sets $\Lambda_n\subset \cF$ of cardinality  $n$ which attain $E_n(K)$. It turns out that the $L_\infty(Y,X)$ norm is difficult to work with and so we replace it by a certain surrogate majorant. In \S \ref{sec:best}, we show that estimating the error of approximating $K$ in the surrogate norm and finding optimal lower sets $\Lambda_n$ in this norm has a simple solution. Namely, given any $\e>0$, the smallest lower set $\Lambda$ for which $\cP_\Lambda$ approximates $K$ to accuracy $\e$ is given by the set of lattice points in a certain simplex $S=S(\e,\rho)$ determined by $\rho$ and $\e$. Therefore, understanding the rate of decay of $E_n(K)$ is equivalent to counting the number of lattice points in these simplices.
  
Of course, counting lattice points in simplices is a well studied problem in number theory where several deep results are known. General results typically only hold for $n$  large when  the number of lattice points can be estimated through the volume of the simplex. In numerical applications, the pre-asymptotic region is the most important  since it corresponds to the only $n$ which can be implemented in computation. Therefore, we focus on counting lattice points when $n$ is small. For results of this type, one needs to have more specific information on the simplices, and therefore on the sequence $\rho$. This leads us to consider specific anisotropic classes that arise in applications. These correspond to sequences $\rho$ which grow polynomially.

For $s>0$, we define the sequence $\rho(s):=(\rho_j(s))_{j\ge1}$ with $\rho_j(s):=(j+1)^s$, $j\ge 1$. The problem of counting lattice points in the simplex associated to this sequence is directly related to  counting the number of multiplicative partitions of integers. One can therefore use the results in \cite{CEP} to do exact and asymptotic analysis for the number of such lattice points. Exact counts on the number of multiplicative partitions of an integer $n$ are known for certain values of $n$. Making such counts becomes numerically more intensive as $n$ increases.
 
It turns out that this situation can be alleviated some by slightly modifying the sequence $\rho(s)$. We modify this sequence to obtain a related sequence $\rho^*(s)$ with the same asymptotic decay as $\rho(s)$. The advantage of this modification is that the number of lattice points of the modified sequence is  related to the number of additive partitions of an integer $n$ rather than the multiplicative partitions. Finding additive partitions is somewhat easier numerically. We show how one can do an exact count of lattice points for $\rho^*(s)$  in \S \ref{sec:pol}.   In \S \ref{sec:asymptotic}, we give an asymptotic analysis for this count by using known asymptotic bounds for the number of additive partitions of integers. In \S \ref{sec:findlambda}, we give some simple recipes for how to find the optimal $\Lambda_n$ for various sequences $\rho$. Finally, in \S \ref{sec:conclusion}, we make some final remarks and compare out results with those in \cite{Tran1}.
  
Let us close this introduction by mentioning that the results in this paper have a large intersection with  several earlier  papers.  As we have already noted, our motivation for the introduction of the spaces $\cB_{\rho,p}$ stems from several works on parametric PDEs, see the survey \cite{CD} and the references therein. Let us also mention \cite{BAS,GO} which study approximation of anisotropic analytic functions in a quite general tensor product framework. There are several papers, most notably \cite{BAS, Tran1, Zech}, that realize that one method to construct approximations for solutions of parametric  PDEs is related to counting lattice points. In~\cite{Tran1}, this count  is done for certain non-simplicial sets as well. We touch more on some of these works later in the paper once our results are formulated.

\section{Anisotropic analyticity and motivation}  
\label{sec:anisoclass}

In this section, we introduce a variety of model classes based on some form of  anisotropic analyticity. We recall that throughout this paper $\rho=(\rho_1,\rho_2,\dots)$ denotes a non-decreasing sequence of positive real numbers with $\rho_1>1$, and $\lim_{j\to\infty}\rho_j=\infty$. We call any sequence with these properties {\it admissible}. 
We recall the Banach  spaces $\ell_\infty(\N)$ of all bounded {\it complex valued} sequences $(z_j)_{j\ge 1}$, with its usual norm $\|z\|_{\ell_\infty{(\N)}}:=\sup_{j\ge 1} |z_j|$.  We let $\cU$ denote the unit ball of $\ell_\infty(\N)$ in this section.

Perhaps the most natural class of anisotropic analytic functions is the following. We start with the complex (open) polydisc $D_\rho$, which consists of all $z=(z_1,z_2,\dots)$, $z_j\in \C$, for which $|z_j|< \rho_j$, $j=1,2,\dots,$ and  define $\overline D_\rho$ as the set of all $z=(z_1,z_2,\dots)$ for which $|z_j|\le  \rho_j$, $z_j\in \C$, $j=1,2,\dots$. We then define
\begin{equation*}
\cH_\rho:=\cH_\rho(X)
\end{equation*}
as the set of all functions $u: \ell_\infty(\mathbb{N}) \rightarrow X$ which are bounded on $\overline  D_\rho$, continuous on $D_\rho$, and holomorphic in each variable $z_j$, $j=1,2, \dots$, on $D_\rho$. We can equip this space with the norm
\bes
\label{Anorm}
{\displaystyle \|u\|_{\cH_\rho}:= \sup_{z\in \overline  D_\rho}\|u(z)\|_X.}
\ees
Because the sequence $\rho$ is non-decreasing, we see that functions in $\cH_\rho$ have more smoothness in the variable $z_j$ as $j$ increases. These spaces are analogous to Hardy spaces.

If $\nu\in\cF$, then the support of $\nu$ is finite and any $u\in \cH_\rho$ has uniquely defined  Taylor coefficients
$$
t_\nu:= \frac{\partial^\nu u(0)}{\nu !}, \quad \nu\in\cF,
$$
where again we are using standard multivariate notation. Note here that the definition of $t_\nu$ requires only the function $u(z_1,\dots,z_N,0,\dots)$, for a suitable finite value of $N$. Since this is an analytic function of a finite number of variables, these coefficients are well defined from the usual theory of functions of a finite number of variables.
 
In what follows, we are interested in representing $u$ in a Taylor series expansion
\be
\label{TE}
u(z)=\sum_{\nu\in\cF} t_\nu z^\nu.
\ee
An important issue is the sense in which the above  Taylor series converges.  For this, we follow Section 3.1 of \cite{CD}. It is shown in that paper that any rearrangement of this   series converges uniformly on $\cU$ whenever  $(\|t_\nu\|_X)_{\nu\in\cF}$ is in $\ell_1(\cF)$. This guarantees that there is a function $v$ defined on $\cU$ such  that any rearrangement of the terms in the series in \eref{TE} converges in $X$ uniformly to $v$. We call this type of convergence {\it uniform unconditional}. 

Convergence of the Taylor series associated to $u$ does not guarantee that its limit is equal to $u$.  For this one requires additional structure. A sufficient condition is that $u$ has the following property:

\vskip .1in
\noindent
{\bf Truncation Property:}
For all $z\in \cU$, we have
\be
\label{tprop}
u(z)=\lim_{N\to\infty}u(z_1,\dots,z_N, 0,\dots).
\ee
\noindent
This property is known to hold for the solutions to parametric PDEs.

Our first observation is that whenever $u$ is   in $\cH_\rho$,  then $u$ has a  bound on its Taylor coefficients.
\begin{lemma}
\label{taylorexpansion}
If $u\in \cH_\rho$,  then 
\begin{enumerate}[label={\rm(\roman*)}]
	\item the Taylor coefficients  $t_\nu\in X$ of $u$ satisfy the bounds  
 \be
 \label{Cauchy}
 \|t_\nu\|_X\le \|u\|_{\cH_\rho} \rho^{-\nu},\quad \nu \in\cF.
 \ee
\item if in addition, $u$ has the {\bf Truncation Property} and $(\|t_\nu\|_X)_{\nu\in\cF}$ is in  $\ell_1(\cF)$, we have
 \bes
 \label{TE1}
  u(z)=\sum_{\nu\in \cF} t_\nu z^\nu,   \quad  \|z\|_{\ell_\infty(\N)}\le 1,
  \ees
  with  uniform unconditional  convergence of the series.  
\end{enumerate}
\end{lemma}
 
\noindent {\bf Proof:} We use a slight modification of  the proof of Lemma 3.14 in \cite{CD} accounting for the fact that the assumptions of the lemma do not guarantee that  $u$ is holomorphic on an open set containing $\overline D_\rho$ as is required in that lemma of \cite{CD}. If we fix $\nu\in \cF$ and any sequence $\delta<\rho$, we claim that
\be
\label{Test}
\|t_\nu\|_X\le \|u\|_{\cH_\rho} \delta^{-\nu}.
\ee
Given $\nu$, let $\{1,\dots,J\}$ contain the support of $\nu$. We consider the  function     $F(z_1,\dots,z_J):=u(\hat z)$, where $\hat z_j:=z_j$ $j=1,\dots,J$, and $\hat z_j$ is zero otherwise. Then $t_\nu$  is the corresponding Taylor coefficient of $F$ and the bound \eref{Test} is derived from Cauchy's formula as in \cite{CD}. Since this bound  holds for any $\delta <\rho$ and the support of $\nu$ is finite, we obtain \eref{Cauchy} by letting $\delta_j\to \rho_j$, for each $j$ in the support of $\nu$. The uniqueness of $t_\nu$ again follows from the fact that $\nu$ has only a finite number of nonzero coordinates and $t_\nu$ is determined by restricting $u$ to the finite number of coordinates corresponding to where $\nu$ is nonzero. This proves (i).  
 
For the proof of (ii), the assumption $(\|t_\nu\|_X)_{\nu\in\cF} \in \ell_1(\cF)$ guarantees the convergence of the Taylor series and then the fact that its sum is $u$ follows easily (see Proposition 2.1.5 in \cite{Zech}). \hfill $\Box$

\vskip .1in

This lemma motivates the definition of the  following class of functions.
\vskip .1in
 
\noindent
{\bf Definition of $\cB_{\rho,\infty}$: } We say that a function $u$ defined on $Y$ and taking values in $X$, is in the space $\cB_{\rho,\infty}:=\cB_{\rho,\infty}(Y)$  if $u$ admits a representation 
\bes
\label{Brep1}
u(y)=\sum_{\nu\in\cF} t_\nu y^\nu,\quad y\in Y,
\ees
with the convergence of the series uniform unconditional on $Y$, and where the $t_\nu=t_\nu(u)\in X$ are unique and satisfy
\bes
\label{defBrho}
\|u\|_{\cB_{\rho,\infty}}:=\sup_{\nu\in\cF} \rho^\nu\|t_\nu\|_X < \infty.
\ees
\vskip .1in
\noindent

Another type of restriction on functions  $u$, derived in the context of parametric PDEs (see \cite{BCM}), is that
\be
\label{BCM}
\sum_{\nu\in \cF}  [\rho^\nu\|t_\nu\|_X ]^2  <\infty.
\ee
This motivates the  general  definition of the following model classes.
  
\vskip .1in
\noindent
{\bf Definition of $\cB_{\rho,p}$:}  For any $0<p\le \infty$, we define the space
$\cB_{\rho,p}$,  as the set  of all $u\in L_\infty(Y,X)$    which admit a representation 
\bes
\label{Brep}
u(y)=\sum_{\nu\in\cF} t_\nu y^\nu,\quad y\in Y,
\ees
with the convergence of the series uniform unconditional on $Y$, and where the $t_\nu=t_\nu(u)\in X$ are unique and satisfy
\be
\label{Brhop}
\|u\|_{\cB_{\rho,p}}:= \left(  \sum_{\nu\in \cF} [\rho^\nu \|t_\nu\|_X]^p\right)^{1/p}= \| (\rho^\nu\|t_\nu\|_X)_{\nu\in\cF}\|_{\ell_p(\cF)} <\infty.
\ee

\vskip .1in
\noindent
Notice that these classes get smaller as $p$ decreases: $\cB_{\rho,p}\subset \cB_{\rho,q}$ when $p\le q$. We study the approximation of the model classes $\cB_{\rho,p}$ in this paper.

We could similarly define anisotropic spaces using other sequence norms in place of $\ell_p$ norms, for example, Lorentz space norms. However, we will not explore this in the present paper.
 
\begin{remark} We have introduced spaces of  anisotropic analytic functions by imposing conditions on Taylor coefficients. One could replace the Taylor basis $y^\nu$, $\nu\in\cF$, by other polynomial bases and define corresponding spaces of analytic functions. A particularly interesting case is when the polynomial basis consists of Legendre polynomials, since such expansions occur naturally in parametric PDEs (see \cite{CCS}).
\end{remark}
 
\section{The approximation of functions in $\cB_{\rho,p}$} 
\label{sec:best}

In this section, we give first estimates for the error in approximating functions in $\cB_{\rho,p}$ by polynomials in $\cP_\Lambda$, with $\Lambda$ a lower set. We follow the ideas in \cite{CM} which treats the case $p=2$. Recall that in this paper we limit our discussion to the approximation of $u$ in the $L_\infty(Y,X)$ norm. This norm is not easy to access especially when $X$ is a general Banach space. However, if $u$ has a Taylor expansion $u(y)=\sum_{\nu\in\cF}t_\nu y^\nu$, $y\in Y$, then it has a simple majorant given by
\bes
\label{majorant}
\|u\|_{L_\infty(Y,X)}\le \sum_{\nu\in\cF}\|t_\nu\|_X=  \sum_{\nu\in\cF}\|t_\nu(u)\|_X =:\|u\|^*.
\ees
The surrogate norm $\|u\|^*$ is defined and finite only if  $u$ has a Taylor expansion valid on $Y$ and $(\|t_\nu\|_X)_{\nu\in\cF}$ is in $\ell_1(\cF)$. We assume that this is the case in going further in this section.  As we shall see below, this assumption is easy to verify when $u\in \cB_{\rho,p}$ under suitable assumptions on $\rho$. This leads us to consider the surrogate error
\be
\label{surrogate}
E_\Lambda^*(u):= \inf_{P\in\cP_\Lambda}\|u-P\|^* =\sum_{\nu\notin\Lambda}\|t_\nu\|_X,
\ee
and similarly
\be
\label{surrogateclass}
E_\Lambda^*(K):= \sup_{u\in K} E^*_\Lambda(u),
\ee
for the surrogate performance on a  compact set $K\subset L_\infty(Y,X)$. Given any  set $\Lambda$, the polynomial 
\begin{equation*}
T_\Lambda(y):=\sum_{\nu\in \Lambda} t_\nu   y^\nu
\end{equation*}
provides an approximation to $u$ which satisfies
\be
\label{eq1}
E_\Lambda(u)\le \|u-T_\Lambda\|_{L_\infty(Y,X)} \le \|u-T_\Lambda\|^*=E_\Lambda^*(u).
\ee

We now describe a simple way to find a lower set from $\cL_n$ which gives the smallest surrogate error for the unit ball $\cU_{\rho,p}$ of $\cB_{\rho,p}$, $0<p\le\infty$, among all lower sets from $\cL_n$. Given the sequence $\rho$ and given any $\e>0$, we define
\be
\label{defLambda} 
\Lambda(\e,\rho):=\{\nu \in\cF: \rho^{-\nu} \ge \e\} =\{\nu \in\cF: \rho^\nu \le \e^{-1}\}.
\ee
Notice that $\Lambda(\e,\rho)$ has the following properties:
\begin{itemize}
\item $\#\Lambda(\e,\rho)<\infty$ whenever $\e>0$, since $\rho$ is non-decreasing, with $\rho_1>1$ and $\lim_{j\to\infty}\rho_j=\infty$;  
\item $\Lambda(\e, \rho)$ is a lower set, since $\mu\leq \nu \Rightarrow \rho^{-\nu}\leq \rho^{-\mu}$;
\item $\Lambda(\e, \rho)\subset\Lambda(\e', \rho)$ whenever $\e'\le\e$.
\end{itemize}
We define  the sequence $(\delta_n)_{n\ge1}=(\delta_n(\rho))_{n\ge1}$ to be a decreasing rearrangement of the sequence $(\rho^{-\nu})_{\nu\in\cF}$. Then, $\#\Lambda(\delta_n,\rho)\ge n$. 
We further define  
\be
\label{defL}
\Lambda_n:=\Lambda_{n,\rho}
\ee 
as any lower set contained in $\Lambda(\delta_n,\rho)$ with cardinality $n$ and which has the property that it contains all $\nu$ for which $\rho^{-\nu}>\delta_n(\rho)$. Such a lower set can be obtained from $\Lambda(\delta_n,\rho)$ by successively removing extreme points and thereby retaining the lower set property. Note that $\Lambda_n$ is not unique because of possible ties in the value of $\rho^{-\nu}$, $\nu\in\cF$.

For any admissible $\rho$ and any $n\ge 1$, we define
\be
\label{defepsq}
  \delta_{n,q}:=\delta_{n,q}(\rho):=
  \begin{cases}
    (\sum_{\nu\notin \Lambda_{n,\rho}} \rho^{-\nu q})^{1/q}=(\sum_{j>n}\delta_j^q)^{1/q}, \quad 
 & \text{if $0<q<\infty$,} \\
    \delta_{n+1}, & q=\infty.
  \end{cases}
\ee
While $\Lambda_n$ need not be unique, we always have a unique value for $\delta_{n,q}(\rho)$ for all choices of $n,q,\rho$.

\begin{theorem}
\label{bestlambdatheorem}
For $0< p\le \infty$,  we have:

\begin{enumerate}[label={\rm(\roman*)}]
\item the set $\Lambda_{n,\rho}$, defined in {\rm \eref{defL}}, minimizes $E^*_\Lambda({ \cU_{\rho,p}})$ over all lower sets $\Lambda\in\cL_n$, and
$$E_n({ \cU_{\rho,p}})\le E_{\Lambda_{n,\rho}}({\cU_{\rho,p}})\le E^*_{\Lambda_{n,\rho}}({ \cU_{\rho,p}});$$
 
\item if $p\ge 1$ and $q$ is the conjugate index to $p$, i.e., $1/p+1/q=1$, then 
$$E_n({\cU_{\rho,p}})\le E_{\Lambda_{n,\rho}}({ \cU_{\rho,p}})\le E^*_{\Lambda_{n,\rho}}({ \cU_{\rho,p}})= \delta_{n,q}.$$
\end{enumerate}
\end{theorem} 

\noindent {\bf Proof:} Let us first make some remarks about the structure of $\cU_{\rho,p}$ that hold for any $0<p\le\infty$ and any admissible $\rho$.  Given any $u\in \cU_{\rho,p}$, we know that $u(y)=\sum_{\nu\in\cF} t_\nu y^\nu$ and the Taylor coefficients $t_\nu$ satisfy
\be
\label{satisfy11}
\|t_\nu\|_X=\alpha_\nu\rho^{-\nu},\ \nu\in\cF,\qquad (\alpha_\nu)_{\nu\in\cF}\in U(\ell_p(\cF)),
\ee 
where $U(\ell_p(\cF))$ is the unit ball of the space $\ell_p(\cF)$. Conversely, let $(\alpha_\nu)_{\nu\in\cF}\in U(\ell_p(\cF))$ be a non-negative sequence and let  $g\in X$ with $\|g\|_X=1$. If we define $t_\nu:= g\rho^{-\nu}\alpha_\nu$, $\nu\in\cF$, then the function $u(y):= \sum_{\nu\in\cF} t_\nu  y^\nu$ will be in $\cU_{\rho,p}$ provided that $(\rho^{-\nu}\alpha_\nu)_{\nu\in\cF}$ is summable.

We first prove (ii) for a fixed $n$ and $p$. We only discuss the case $p>1$. The case $p=1$ is proved in a similar way. The two inequalities in (ii) are obvious from the definitions \eqref{classerrorlower},~\eqref{surrogateclass}, and~\eqref{eq1}, and so we only need to show that $E^*_{\Lambda_{n,\rho}}({ \cU_{\rho,p}})=\delta_{n,q}$. Let  $u\in \cU_{\rho,p} $ with $u(y)=\sum_{\nu\in \cF} t_\nu y^\nu$. It follows from \eref{surrogate} with $T_{\Lambda_{n,\rho}}(y):=\sum_{\nu\in \Lambda_{n,\rho}}t_\nu y^\nu$ and H\"older's inequality that
\be
\label{bt1}
E^*_{\Lambda_{n,\rho}}(u) = \|u-T_{\Lambda_{n,\rho}}\|^*=  \sum_{\nu\notin \Lambda_{n,\rho}} \|t_\nu\|_X =\sum_{\nu\notin \Lambda_{n,\rho}}\|t_\nu\|_X\rho^\nu \rho^{-\nu}\le \|u\|_{\cB_{\rho,p}} \delta_{n,q}\le \delta_{n,q}.
\ee
To prove that $E_{\Lambda_{n,\rho}}^*(\cU_{\rho,p})\ge \delta_{n,q}$, we construct a function $\tilde u \in \cU_{\rho,p}$ for which $E_{\Lambda_{n,\rho}}^*(\tilde u)=\delta_{n,q}$. First assume that $\delta_{n,q}$ is finite, so that there is a nonnegative  sequence  $(c_\nu)_{\nu\in\cF}$ in the unit ball of $\ell_p(\cF)$ for which $\sum_{\nu\notin \Lambda_{n,\rho}} c_\nu\rho^{-\nu} =\delta_{n,q}$. Then, as in our lead remarks, we let $g\in X$ with $\|g\|_X = 1$ and define $t_\nu:= c_\nu \rho^{-\nu}g$. Then, we have  
$$\tilde u(y):=\sum_{\nu\in\cF} t_\nu y^\nu,\quad y\in Y,$$ 
is in $\cU_{\rho,p}$. Note here we use the fact that $(\|t_\nu\|_X)_{\nu\in\cF}$ is in $\ell_1(\cF)$. Since $E_{\Lambda_{n,\rho}}^*(\tilde u)=\delta_{n,q}$, we have finished the proof of (ii) in the case that $\delta_{n,q}$ is finite. If $\delta_{n,q}=\infty$, the same argument as above shows that there is a $\tilde u$ for which $E_{\Lambda_{n,\rho}}^*(\tilde u)$ is as large as we wish. Therefore (ii) holds in this case as well.

Now, consider the proof of (i). The inequalities stated in (i) are all obvious and so we need only show that $\Lambda_{n,\rho}$ minimizes $E_\Lambda^*(\cU_{\rho,p})$ over all lower set $\Lambda\in\cL_n$. To prove this, we first consider the case $p\le 1$. If $\Lambda\in\cL_n$, then by our lead remarks
\be
\label{classerror2}
E^*_\Lambda(\cU_{\rho,p})=\sup_{u\in \cU_{\rho,p}} \sum_{\nu\notin \Lambda}\|t_\nu(u)\|_X= \sup_{\alpha\in\ U(\ell_p(\cF))}\sum_{\nu\notin \Lambda} \alpha_\nu \rho^{-\nu}.
\ee
Here, we use the fact that $(\rho^{-\nu}\alpha_\nu)_{\nu\in\cF}$ is summable because $(\alpha_\nu)_{\nu\in\cF}$ is in $\ell_1(\cF)$. The minimum of \eref{classerror2} over all $\Lambda$ is achieved by taking $\Lambda=\Lambda_{n,\rho}$.   
 
Finally, we have to prove (i) in the case $1<p\le \infty$. Let us first recall that there is an enumeration $\nu(n)$, $n\ge 1$, of all of the $\nu\in\cF$, such that $\delta_n=\rho^{-\nu(n)}$, $n\ge 1$, and such that $\Lambda_{n,\rho}=\{\nu(1),\dots,\nu(n)\}$. We suppose $\Lambda$ is any lower set with $\#\Lambda=n$. From the definition of $\Lambda_{n,\rho}$, we have
\be
\label{bt2}
\sum_{\nu\notin \Lambda} \rho^{-\nu q}\ge \delta_{n,q}^q.
\ee
Using the same construction as in the proof of (ii), we can find $\tilde u \in \cU_{\rho,p}$ with
$$E_\Lambda^*(\tilde u)^q = \sum_{\nu\notin \Lambda} \rho^{-\nu q}\ge \delta_{n,q}^q,$$
which thereby proves (i).\hfill $\Box$

\begin{cor}
\label{firsttcorollary}  
For $1\le p\le \infty$, let $q$ be the conjugate index to $p$. Then whenever $\delta_{n,r}$ is finite for some $0<r < q$, we have
\be
\label{fc1}
E_n(\cU_{\rho,p})\le E_{\Lambda_{n,\rho}}(\cU_{\rho,p})\le E_{\Lambda_{n,\rho}}^*(\cU_{\rho,p})\le \delta_{n+1}^{1-r/q} \delta_{n,r}^{r/q},\quad n\ge 1.
\ee
In particular, we have
\be
\label{fc2}
E_{\Lambda_{n,\rho}}^*(\cU_{\rho,p})\le (n+1)^{-1/r+1/q}\|(\rho^{-\nu})_{\nu\in\cF}\|_{\ell_r(\cF)},\quad n\ge 1.
\ee
\end{cor}

\noindent {\bf Proof:} The case where $p=1$ (resp. $q=\infty$) follows from (ii) of Theorem \ref{bestlambdatheorem}, since $\delta_{n,\infty}=\delta_{n+1}$. So we can assume $p > 1$ and $q < \infty$. Since the sequence $(\delta_n)_{n\ge1}$ is non-increasing, we have
\be
\label{dec}
\delta_{n,q}^q=\sum_{j>n}\delta_j^q\le  \delta_{n+1}^{q-r}\sum_{j>n} \delta_j^r\le \delta_{n+1}^{q-r}\delta_{n,r}^r.
\ee
Because of (ii) in Theorem \ref{bestlambdatheorem}, taking a $q$-th root proves \eref{fc1}. To show  \eref{fc2}, we
use the fact that  $\delta_{n,r}\le \|(\rho^{-\nu})_{\nu\in\cF}\|_{\ell_r(\cF)}$, along with the standard estimate
$$(n+1)\delta_{n+1}^r\le \sum_{j=1}^{n+1} \delta_j^r\le  \|(\rho^{-\nu})_{\nu\in\cF}\|_{\ell_r(\cF)}^r.$$
Inserting these into \eref{dec}, we obtain
$$
\delta_{n,q}^q\le \delta_{n+1}^{q-r}\delta_{n,r}^r\le (n+1)^{-\frac{q-r}{r} } \|(\rho^{-\nu})_{\nu\in\cF}\|^{q}_{\ell_r(\cF)}=
(n+1)^{-\frac{q}{r}+1} \|(\rho^{-\nu})_{\nu\in\cF}\|^{q}_{\ell_r(\cF)},
$$
and the proof is complete.
\hfill $\Box$

\begin{remark}
\label{finited}
We can define the above space $\cB_{\rho,p}$ also in the case $\rho=(\rho_1,\dots,\rho_d)$ with $d$ finite. The results of this section hold equally well in this case.    
 \end{remark}

\section{The sequence $\delta_n(\rho) $} 
\label{sec:gen}
 
First, let us observe that in order for the $\delta_{n,q}$ from~\eref{defepsq} to be finite, and therefore Theorem \ref{bestlambdatheorem} to be meaningful, we need that the sequence $(\rho^{-\nu})_{\nu \in \cF} \in \ell_q(\cF)$, which is the same as asking that  $(\delta_n(\rho))_{n\geq 1} \in \ell_q(\N)$. The following lemma shows that this is the case if and only if $(\rho_j^{-1})_{j\ge1} \in\ell_q(\N)$.

\begin{lemma}
\label{lqlemma}
Let $0<q\le  \infty$. Then the sequence $(\rho^{-\nu})_{\nu\in\cF} \in\ell_q(\cF)$ if and only if the sequence $(\rho_j^{-1})_{j\ge1}{ \in}\ell_q(\mathbb{N})$. Moreover, the two norms are related in the following way:
 
\begin{enumerate}[label={\rm(\roman*)}]
\item when $q=\infty$, we have ${\|(\rho^{-\nu})_{\nu\in\cF}\|}_{\ell_\infty{ (\cF)}}=1,$
\item when  $0<q<\infty$, we have
\bes
e^{{ \|(\rho_j^{-1})_{j\geq1}\|}^q_{\ell_q{ (\N)}}}\le { \|(\rho^{-\nu})_{\nu\in\cF}\|}^q_{\ell_q{ (\cF)}}\le e^{(1-\rho_1^{-q})^{-1}{ \|(\rho_j^{-1})_{j\geq1}\|}^q_{\ell_q{ (\N)}}}.
\ees
\end{enumerate}
\end{lemma}

\noindent {\bf Proof:} The case $q=\infty$ is trivial. When $q<\infty$, we have
\be
\label{inlp} 
\| (\rho^{-\nu})_{\nu\in\cF}\|_{\ell_q{ (\cF)}}^q=\sum_{\nu\in\cF}\rho^{-q\nu}= \prod_{j\ge 1}\left( \sum_{k=0}^\infty \rho_j^{-kq}\right) =\prod_{j\ge 1}(1-\rho_j^{-q})^{-1}.
\ee
Taking logarithms, we have from the mean value theorem that
\bes
\label{prod}
\ln (\prod_{j=1}^\infty (1-\rho_j^{-q})^{-1}) = - \sum_{j=1}^\infty \ln (1-\rho_j^{-q})=\sum_{j=1}^\infty (1-\xi_j)^{-1} \rho_j^{-q},
\ees
where the $\xi_j\in (0,\rho_j^{-q})\subset (0,\rho_1^{-q})$, $j=1,2,\dots$. Since $1< (1-\xi_j)^{-1} < (1-\rho_1^{-q})^{-1}$, it follows that
\bes
\label{comparenorms}
e^{{ \|(\rho_j^{-1})_{j\geq1}\|}^q_{\ell_q{ (\N)}}}\le { \| (\rho^{-\nu})_{\nu\in\cF}\|}^q_{\ell_q{ (\cF)}}\le e^{(1-\rho_1^{-q})^{-1}{ \|(\rho_j^{-1})_{j\geq1}\|}^q_{\ell_q{ (\N)}}}.
\ees
This proves item (ii) in the lemma, and likewise shows that the product in~\eqref{inlp} converges if and only if  $(\rho_j^{-1})_{j\ge1}{ \in}\ell_q{ (\N)}$. \hfill $\Box$

\begin{remark}
\label{remarkcm}
The upper bound established in the above lemma can be found in \cite{CM}.
\end{remark}

The error estimates derived in \S\ref{sec:best} for approximation by polynomials on lower sets depend crucially on the sequence $(\delta_n(\rho))_{n\ge1}$, and are achieved by choosing the lower set $\Lambda_n=\Lambda_{n,\rho}$.  This leads to two central issues:  
\begin{enumerate}[label={\rm(\roman*)}]
\item establishing sharp a priori estimates for $\delta_n(\rho)$ given the sequence $\rho$;
\item efficient algorithms for generating the sets $\Lambda_n$.
\end{enumerate}
We discuss item (ii) in \S \ref{sec:findlambda}, and here we discuss first item (i). We begin this section with methods for bounding $(\delta_n(\rho))_{n\ge1}$ which hold for any admissible sequence $\rho$.   
\begin{remark}
\label{inverseremark}
In order to compute $\delta_n(\rho)$ or its asymptotic decay as $n\to\infty$, we study $\#\Lambda(\e,\rho)$, $0<\e\le1$. This function of $\e$ takes integer values and increases as $\e$ goes to zero. Hence, it is a piecewise constant function and  $(\delta_n(\rho))_{n\ge1}$ is the decreasing  sequence of the breakpoints $\e_1, \e_2,\ldots,$ of $ \#\Lambda(\e,\rho)$, where each value $\e_i$ is repeated $\#\Lambda(\e_{i+1},\rho)-\#\Lambda(\e_i,\rho)$ times and $\delta_1(\rho)=\ldots=\delta_{k_1}(\rho)=1$ with $k_1=\#\Lambda(1,\rho)$. 
\end{remark}

Since $\displaystyle{\lim_{j\to\infty} \rho_j=+\infty}$, there is a $D=D(\e)$ such that $\rho_j^{-1}<\e$, $j>D$.  It follows that any $\nu\in \Lambda(\e,\rho)$ has support in $\{1,2,\dots,D\}$. Moreover, if we write $\e=e^{-M}$, then taking logarithms we see that $\nu\in \Lambda(\e,\rho)$ if and only if $\nu$ satisfies
\bes
\label{iffbig}
\sum_{j=1}^D  \nu_j  \frac{ \ln \rho_j}{M}\le 1.
\ees
Hence, $\nu\in\Lambda(\e,\rho)$ if and only if $\nu$ is supported on $\{1,2,\dots,D\}$, and $(\nu_1,\dots,\nu_D)$ is a lattice point in the simplex
\be
\label{simplex}
S:=S(a_1,\dots,a_D):=\{(x_1,\dots,x_D): \ x_i\ge 0, \ i=1,\dots,D, \ {\rm and} \ \sum_{j=1}^D\frac{ x_j}{a_j}\le 1\}, 
\ee
where
$$a_j:=  \frac{M}{\ln \rho_j}\ge 1, \quad j=1,\dots,D.$$
Estimating the number of lattice points in such a simplex is a classical problem in number theory and combinatorics. Let us first note that the volume (measure) of $S$ is
\be
\label{volume}
{\rm vol}(S)=|S|= \frac{\prod_{j=1}^D a_j}{D!}.
\ee
We recall the following general upper bound (see \cite{B-D, YZ}) for the number $\#\Lambda (S)$ of $\nu\in\mathbb{N}_0^D$ such that $\nu\in S$:
\be
\label{upperbound}
\frac{\prod_{j=1}^D a_j}{D!}\le  \#\Lambda(S)\le (1+a)^D\frac{\prod_{j=1}^D a_j}{D!}, \quad a:=\sum_{j=1}^D a_j^{-1}.
\ee
Note that the right side of \eref{upperbound} is inflated by a factor of ($1+a)^D$ when compared with the volume of $S$. We use this result to prove the following lemma.
  
\begin{lemma}
\label{countlemma}
Let $\rho$ be any admissible sequence. Given $\e=e^{-M}$, where $M>0$, let $D$ be the last integer $j$ for which $\rho_j\le e^M$. Then, for the set $\Lambda(\e,\rho)$ of all  $\nu\in\cF$ such that $\rho^{-\nu}\ge \e$, we have
\be
\label{boundbig}
\#\Lambda(\e,\rho)\le \frac{ (M+L)^D }{D!} [\prod_{j=1}^D \ln \rho_j]^{-1},\quad \hbox{where}\quad L:=L(\rho):=\sum_{j=1}^D  \ln \rho_j .
\ee
\end{lemma}
 
\noindent{\bf Proof:} From \eref{upperbound} with $a_j=\frac{M}{ \ln \rho_j}$, $j=1,\dots,D$, we have
\bes
\label{count1}
\#\Lambda(\e,\rho)\le   \frac{M^D} {D!}\left (1+\frac{L}{M}\right)^D[ \prod_{j=1}^D \ln \rho_j] ^{-1},
\ees
which is equivalent to \eref{boundbig}.   \hfill $\Box$
 
Let us make some remarks that will clarify when the bound in the lemma is effective and when it is  deficient.  
First of all, if $d$  is finite and the sequence $(\rho_j)_{j= 1}^d$ is  fixed, then the set $ \Lambda(\e,\rho)$, $\e=e^{-M}$, is the set of lattice points $\N^d_{0}/M$ in the fixed simplex $S^*:=S(1/\ln \rho_1,\dots, 1/\ln \rho_d)$. If we let $M$ tend to infinity (which corresponds to $\e\to 0$), we see that $D=d$ provided $M$ is large enough, and $\#\Lambda(e^{-M},\rho )$ behaves like $M^d$ times the measures of $S^*$. This is in agreement with the bound \eref{boundbig} because the inflation factor $(1+L/M)^D=(1+L/M)^d$ tends to one as $M\to\infty$. So this bound is good for finite $d$, provided the error we seek is small. However, there is a transition before this asymptotic kicks in where the upper bound provided by the lemma is not effective.
 
To see this, we consider one example which is central to this paper. We consider the sequence $\rho:=(j+1)_{j=1}^d$, with  $d$ finite. We take as our target error $\e:=1/(d+1)$, i.e. $M= \ln (d+1)$. Then $D=d$ and the upper bound for $\#\Lambda(\e,\rho)$ provided by Lemma \ref{countlemma} is
\be
\label{lambdaestimate}
\frac{(\ln (d+1) +\ln(d+1)!)^d}{d!   \prod_{j=1}^d \ln(j+1)}=:B(d), 
\ee
where we used the fact that 
\bes
\label{estimateL}
L= \sum_{j=1}^d  \ln (j+1)=  \ln (d+1)!, \quad \prod_{j=1}^d \ln \rho_j =   \prod_{j=1}^d \ln (j+1).
\ees
Since $\ln(x)$ is a concave function, we have 
$$
\ln(d+1)!=\sum_{j=2}^{d+1}\ln j\ge \int_1^{d+1} \ln t\,dt  \ge \frac{d \ln (d+1)}{2}.
$$
Therefore, we have 
\begin{eqnarray}
\nonumber
B(d)&\ge& \frac{[\ln(d+1)!]^d}{d!\prod_{j=1}^d \ln(j+1)}\ge 
\frac{d^d[\ln(d+1)]^d}{2^d d!\prod_{j=1}^d \ln(j+1)}=\frac{d^d}{2^dd!}\prod_{j=1}^d \frac{\ln(d+1)}{\ln(j+1)}
\\ \nonumber
&\ge &\frac{d^d}{2^dd!}\ge \frac{1}{e\sqrt{d}}\left (\frac{e}{2}\right)^d,
\nonumber
\end{eqnarray}
where we used Stirling's formula. Thus, if we want an error $\e=1/(d+1)$ in this particular example, the best bound that Lemma \ref{countlemma} can provide for the size of $\Lambda(\e,\rho)$ is exponential in $d$. In contrast, in Lemma \ref{aslemma} from the following section, we give a much more favorable bound.

\section{Analysis of $\delta_n(\rho) $ when $\rho$ has polynomial growth}
\label{sec:pol}

As we have just observed, the bounds of the previous section for $\delta_n(\rho)$ are generally far from sharp. We can establish sharper bounds, and even compute $\delta_n(\rho)$ exactly, if we have more information on the sequence $\rho$. In this section, we give such an analysis when the sequence $\rho$ has polynomial growth.
 
Recall that for $s>0$, we defined the sequence $\rho(s):=((j+1)^s)_{j\ge 1}$. In some parts of our analysis, it is useful to slightly modify this sequence. Accordingly, we introduce the following modified sequence  $\rho^*(s)$, $s>0$, defined as follows. If $I_1:=\{1,2\}$ and $I_k:=\{j:\ 2^{k-1}<j\le 2^k\}$, $k\ge 2$, then
\be \label{rhos}
\rho_j^*(s):=2^{ks},  \quad  j\in I_k,\quad k=1,2,\dots.
\ee
Note that the sequence $\rho_j^*(s)$ increases like $j^s$. Moreover, $\#I_1=2$ and $\#I_k=2^{k-1}$ for $k\ge 2$.
  
Given any $\e$, we want to determine the cardinality of the set $\Lambda(\e,\rho(s))$ or its counterpart $\Lambda(\e, \rho^*(s))$, i.e., how many $\nu$ satisfy the inequality $[\rho ^*(s)] ^{-\nu} \ge  \e$. According to Remark \ref{inverseremark}, the decay rate of $\delta_n(\rho^*(s))$ can then be derived from this knowledge. Let us note that for these two sequences, we have
\be
\label{sone}
\Lambda(\e^s,\rho(s))=\Lambda(\e,\rho(1)),\quad  \Lambda(\e^s,\rho^*(s))=\Lambda(\e,\rho^*(1)),\quad s>0,
\ee
and so it is enough to analyze the case $s=1$. We therefore take $s=1$ in the estimates on cardinality that follow.
 
As $\e$ decreases, the cardinality of $\Lambda(\e,\rho(1))$ increases. While it is interesting to understand how this cardinality grows asymptotically when $\e$ tends to zero, in numerical scenarios it is important to keep this cardinality small. 

\subsection{Exact formulas for $ \#\Lambda(\e,\rho^*(1))$}

Exact formulas  for the cardinality of $\Lambda(\e,\rho(1))$ can be  given in terms of the multiplicative partitions of natural numbers (see \cite{CEP} and Remark 3.18 in \cite{CD}). In theory, these formulas  allow  the precise computation of $\#\Lambda (\e,\rho(1))$ provided that this cardinality is not too large. However, this computation is very intense and in fact, to our knowledge, has not been done. It turns out that these computations are simpler if one uses the sequence $\rho^*(1)$ instead of $\rho(1)$. This stems from the fact that $\rho^*(1)^\nu$ is always an integer power of two.  For this reason, we focus on this sequence for the remainder of this section. We begin by showing how one can do an exact count of the multiindices in the simplex associated to $\rho^*(1)$.
 
For any $m\in \N_0$, we define  
\be
\label{sumn}
S_m:=\{\nu\in\cF:\  \rho^*(1)^\nu= 2^m\}.
\ee
The set $S_0$ contains only the zero sequence and hence $\#S_0=1$. We want to determine the cardinality of the sets  $S_m$, $m\ge 1$. This is the same as finding how many $\nu\in \cF$ satisfy \eref{defLambda}, since if we denote by 
$$m(\e):=\left\lfloor\log_2\left (\frac{1}{\e}\right) \right\rfloor,$$
we have that 
\be
\label{pq}
\#\Lambda(\e,\rho^*(1))= \#\left( \bigcup_{k=0}^{m(\e)} \{ \nu \in \cF ~:~  \rho^*(1)^\nu= 2^k  \}\right)=\sum_{k=0}^{m(\e)}\#S_k.
\ee

Let us first note that if $\nu$ has a nonzero component $\nu_j>0$ for some $j> 2^m$, then $\rho^*(1)^\nu>2^m$ and so $\nu$ is not in $S_m$. Hence, any $\nu\in S_m $ is supported on $\{1,\dots, 2^m\}$. We decompose the set $\{1,\dots,2^m\}=\bigcup_{k=1}^{m} I_k$, and given any $\nu $, we define 
$$N_k(\nu):=\sum_{j\in I_k}\nu_j,$$
which we think of as the energy of $\nu$ on $I_k$. Therefore, for any $\nu\in\ S_m$, we have
\be
\label{know11}
\quad \sum_{k=1}^m  k N_k(\nu) = m.
\ee
Note that there are only certain sequences $(N_1,\dots,N_m)$ which satisfy \eref{know11}. We denote the collection of all such sequences by $\cQ_m$,
$$\cQ_m:=\{(N_1,\ldots,N_m):\,\,\sum_{k=1}^m  k N_k= m,\,\,N_i\in\N_0,  i=1,\ldots,m\}.$$
The sequences in $\cQ_m$ are related to the additive partitions of $m$, which are decompositions of $m\in\mathbb{N}$ into $m=m_1+\dots+m_j$, where the $m_j\in \N$ and where the order of the appearance of an $m_j$ does not matter. 

There is a one to one correspondence between the elements in $\cQ_m$ and  additive partitions of $m$. Indeed, any  additive partition  $(m_1, m_2,\dots,m_j)$ of $m$ corresponds to a sequence \\ $(N_1,\dots, N_k,\dots,N_m)\in\cQ_m$, 
where $N_k$ is the number of appearances of $k$ in $(m_1,\dots,m_j)$. Conversely, any $(N_1,\dots,N_m)$ for which $\sum_{k=1}^m kN_k=m$ corresponds to the unique additive partition of $m$, where $1$ appears $N_1$ times, $2$ appears $N_2$ times and so on. Thus $q(m):=\#\cQ_m$ is the additive partition number of $m$.

The following theorem gives an exact count for the cardinality of $S_m$, and hence the cardinality of the set $\Lambda(\e,\rho^*(1))$.      
\begin{theorem}
\label{cardtheorem}
For $m\ge 1$, the cardinality of $S_m$ is given by
\be
\label{countSm}
\#S_m = \sum_{(N_1,\dots,N_m)\in \cQ_m} \prod_{k=1}^m  \binom {N_k-1+\#I_k}{N_k}.
\ee
Moreover, for every $\e>0$, 
\be
\label{exact}
\#\Lambda(\e,\rho^*(1)) = \sum_{k=0}^{m(\e)} \#S_k=1+\sum_{k=1}^{ m(\e)}
\sum_{(N_1,\dots,N_k)\in \cQ_{k}}\prod_{j=1}^k \binom {N_j-1+\#I_j}{N_j},
\ee
where $m(\e):=\lfloor\log_2\left (\frac{1}{\e}\right) \rfloor.$
\end{theorem}

\noindent {\bf Proof:} For any fixed $(N_1,\dots,N_m)\in \cQ_m$, we define
\bes
\label{defineLambda}
\Gamma(N_1,\dots,N_m):=\{ \nu\in S_m: N_k(\nu)=N_k,\ k=1,\dots,m\}.
\ees
Now, for each $k=1,\ldots,m$, we count all possible $\nu$ satisfying
$$N_k=N_k(\nu)=\sum_{j\in I_k}\nu_j.$$
Since $\nu_j\in\N_0$, the latter cardinality can be viewed as the  number of ways one can place $N_k$ indistinguishable balls into $\#I_k$ distinguishable boxes so that some boxes can remain empty. The answer to this combinatorial problem is known to be $\binom {N_k-1+\#I_k}{N_k}$ (see \cite{R}). Therefore, the cardinality of $\Gamma(N_1,\dots,N_m)$ is the product of these binomial coefficients:
\be
\label{know21}
\#\Gamma(N_1,\dots,N_m)=  \prod_{k=1}^m \binom {N_k-1+\#I_k}{N_k}.
\ee
Equation \eref{countSm} now follows from the definitions of $S_m$ and $\cQ_m$ and \eref{know21}. The last statement in the theorem follows from \eref{pq} and \eref{countSm}.
\hfill $\Box$

Theorem~\ref{cardtheorem} gives an exact formula for $\#\Lambda(\e,\rho^*(1))$ for  any $\e$, since
$$\Lambda(\e,\rho^*(1))=\Lambda(2^{-m},\rho^*(1)) \quad \text{for}\quad \e\in(2^{-(m+1)},2^{-m}].$$
Note that the sequence $(\delta_n(\rho^*(1)))_{n\ge1}$ is then given by $\delta_1(\rho^*(1))=1$ and, for $m=1,2,\ldots$, 
\be
\label{dt}
\delta_n(\rho^*(1))=2^{-{m}}\quad \text{for}\quad n=\#\Lambda(2^{-m+1},\rho^*(1))+1,\ldots,
\#\Lambda(2^{-m+1},\rho^*(1))+\#S_{m}.
\ee

Moreover, since $\Lambda(\e^s,\rho^*(s))=\Lambda(\e,\rho^*(1))$, $s>0$, we similarly derive that
\be
\label{lambdas}
\Lambda(\e,\rho^*(s))=\Lambda(2^{-m},\rho^*(1)) \quad \text{for}\quad \e\in(2^{-(m+1)s},2^{-ms}],
\ee
and $(\delta_n(\rho^*(s)))_{n\ge1}$ is then given by $\delta_1(\rho^*(s))=1$ and, for
$m=1,2,\ldots$, 
\be
\label{dts}
\delta_n(\rho^*(s))=2^{-{ms}}\quad \text{for}\quad n=\#\Lambda(2^{-m+1},\rho^*(1))+1,\ldots,
\#\Lambda(2^{-m+1},\rho^*(1))+\#S_{m}.
\ee

In Table {\rm \ref{Lambdatb}}, we present the computed cardinality $\#\Lambda(2^{-ms},\rho^*(s))$ for values of $m$ in the range $0\le m\le 10$ and $s=1,2,3,4$. 

\begin{table}[h!]
\centering
\begin{tabular}{| c | c | c | c | c | c |}
	\hline
	\multirow{2}{*}{$ m $}   &  \multirow{2}{*}{$ \#\Lambda(2^{-ms},\rho^*(s)) $} &   \multicolumn{4}{c|}{$2^{-ms}$} \\
	\cline{3-6}
	& & $s=1$ & $s=2$ & $s=3$ & $s=4$\\
	\hline
	$ 0 $ & $ 1 $ & $ 1 $ & $ 1 $ & $ 1 $ & $ 1 $ \\
	$ 1 $ & $ 3 $ & $  5.0000\times 10^{-1} $ & $  2.5000\times 10^{-1} $ & $  1.2500\times 10^{-1} $ & $  6.2500\times 10^{-2} $ \\
	$ 2 $ & $ 8 $ & $  2.5000\times 10^{-1} $ & $  6.2500\times 10^{-2} $ & $  1.5625\times 10^{-2} $ & $  3.9062\times 10^{-3} $ \\
	$ 3 $ & $ 20 $ & $  1.2500\times 10^{-1} $ & $  1.5625\times 10^{-2} $ & $  1.9531\times 10^{-3} $ & $  2.4414\times 10^{-4} $ \\
	$ 4 $ & $ 50 $ & $  6.2500\times 10^{-2} $ & $  3.9062\times 10^{-3} $ & $  2.4414\times 10^{-4} $ & $  1.5259\times 10^{-5} $ \\
	$ 5 $ & $ 122 $ & $  3.1250\times 10^{-2} $ & $  9.7656\times 10^{-4} $ & $  3.0518\times 10^{-5} $ & $  9.5367\times 10^{-7} $ \\
	$ 6 $ & $ 298 $ & $  1.5625\times 10^{-2} $ & $  2.4414\times 10^{-4} $ & $  3.8147\times 10^{-6} $ & $  5.9605\times 10^{-8} $ \\
	$ 7 $ & $ 718 $ & $  7.8125\times 10^{-3} $ & $  6.1035\times 10^{-5} $ & $  4.7684\times 10^{-7} $ & $  3.7253\times 10^{-9} $ \\
	$ 8 $ & $ 1723 $ & $  3.9062\times 10^{-3} $ & $  1.5259\times 10^{-5} $ & $  5.9605\times 10^{-8} $ & $  2.3283\times 10^{-10} $ \\
	$ 9 $ & $ 4101 $ & $  1.9531\times 10^{-3} $ & $  3.8147\times 10^{-6} $ & $  7.4506\times 10^{-9} $ & $  1.4552\times 10^{-11} $ \\
	$ 10 $ & $ 9712 $ & $  9.7656\times 10^{-4} $ & $  9.5367\times 10^{-7} $ & $  9.3132\times 10^{-10} $ & $  9.0949\times 10^{-13} $ \\
	\hline
\end{tabular}
\caption{Computed cardinality of $\Lambda(2^{-ms},\rho^*(s))$ using Theorem \ref{cardtheorem} and the relation (\ref{lambdas}). When $\rho =\rho^*(s)$, the table gives the cardinality of the lower set needed to achieve accuracy $2^{-ms}$. Refer to Remark~\ref{rmk} to deduce estimates on the errors $E_n(\cU_{\rho^*(s),1})$.} \label{Lambdatb}
\end{table} 

\begin{remark}  
\label{rmk}
If we combine this theorem with Theorem {\rm \ref{bestlambdatheorem}} and \eref{dts}, we determine the optimal error and best lower set for approximating any of the spaces $\cB_{\rho^*(s),p}$, provided the error is measured in the surrogate norm rather than the true $L_\infty(Y,X) $ norm. Of course, it gives an upper bound on the performance in the $L_\infty(Y,X)$ norm, that is for $n\ge 1$ we have
$$E_n(\cU_{\rho^*(s),p})\le \left(\sum_{j>n}\delta_j^q\right)^{1/q}, \quad \frac{1}{q}+\frac{1}{p}=1, \quad 1\le q<\infty,$$
and
$$E_n(\cU_{\rho^*(s),1})\le \delta_{n+1},$$
where the sequence $(\delta_n)_{n\ge1}=(\delta_n(\rho^*(s)))_{n\ge1}$ is given by \eref{dts}. The efficiency of the algorithm is determined by the cardinality of $\Lambda(2^{-ms},\rho^*(s)) = \Lambda(2^{-m},\rho^*(1))$, given in Table \ref{Lambdatb}. In particular, let us suppose the user desires to approximate a function in $\cU_{\rho^*(s),1}$ with accuracy $10^{-3}$. Because $\delta_{n+1}=2^{-ms}$ for $n=\#\Lambda(2^{-m+1},\rho^*(1))$ according to \eqref{dts}, when $s=1$, we need $m=10$ and thus a set $\Lambda$ of cardinality 4101 achieves this accuracy. Similarly, a sufficient cardinality for $\Lambda$ is 50 when $s=2$; 20 for $s=3$;  8 for $s=4$.
\end{remark}

\begin{figure}[htbp]
	\centering
	\includegraphics[width=0.45\textwidth]{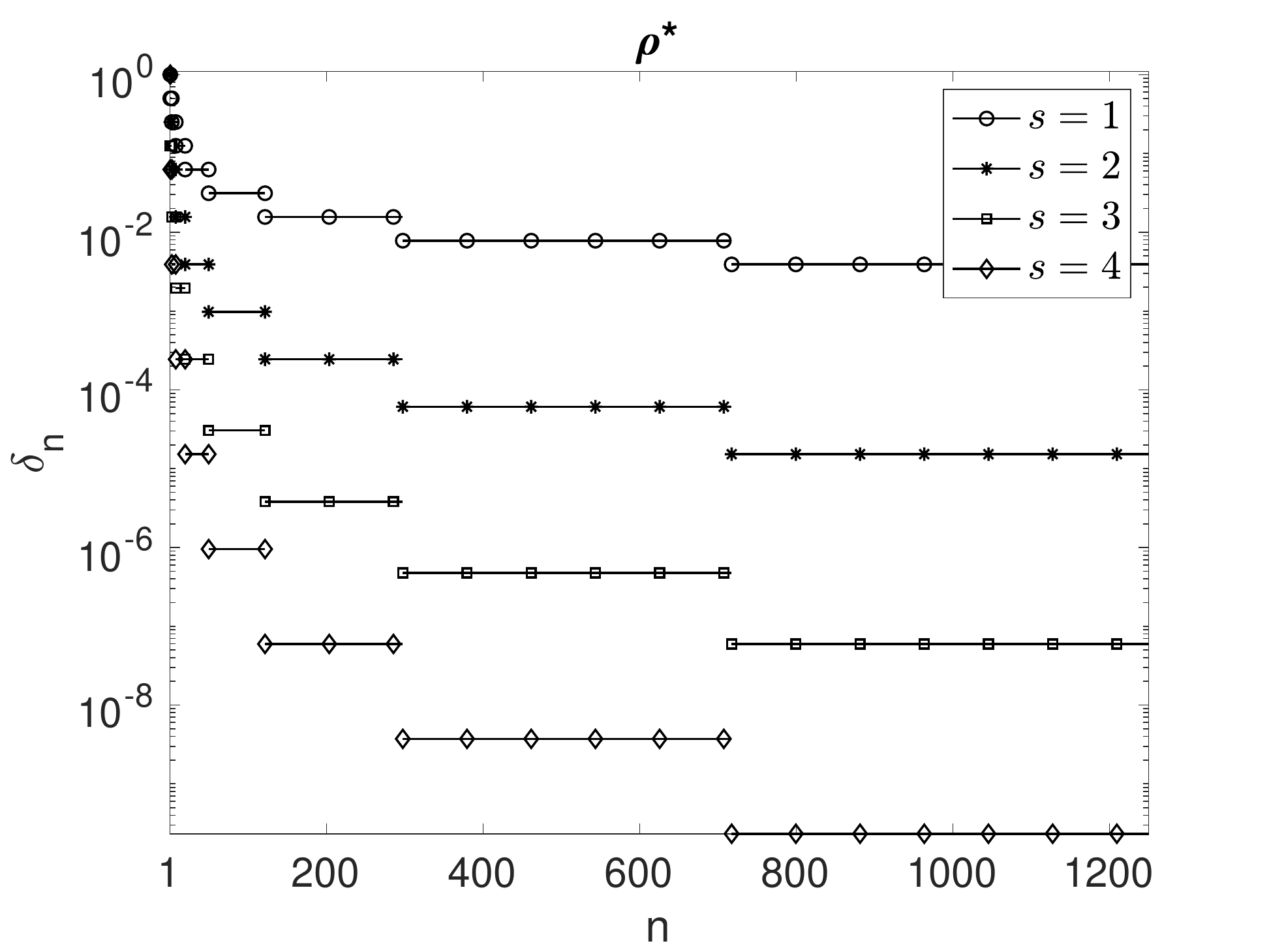}
	\includegraphics[width=0.45\textwidth]{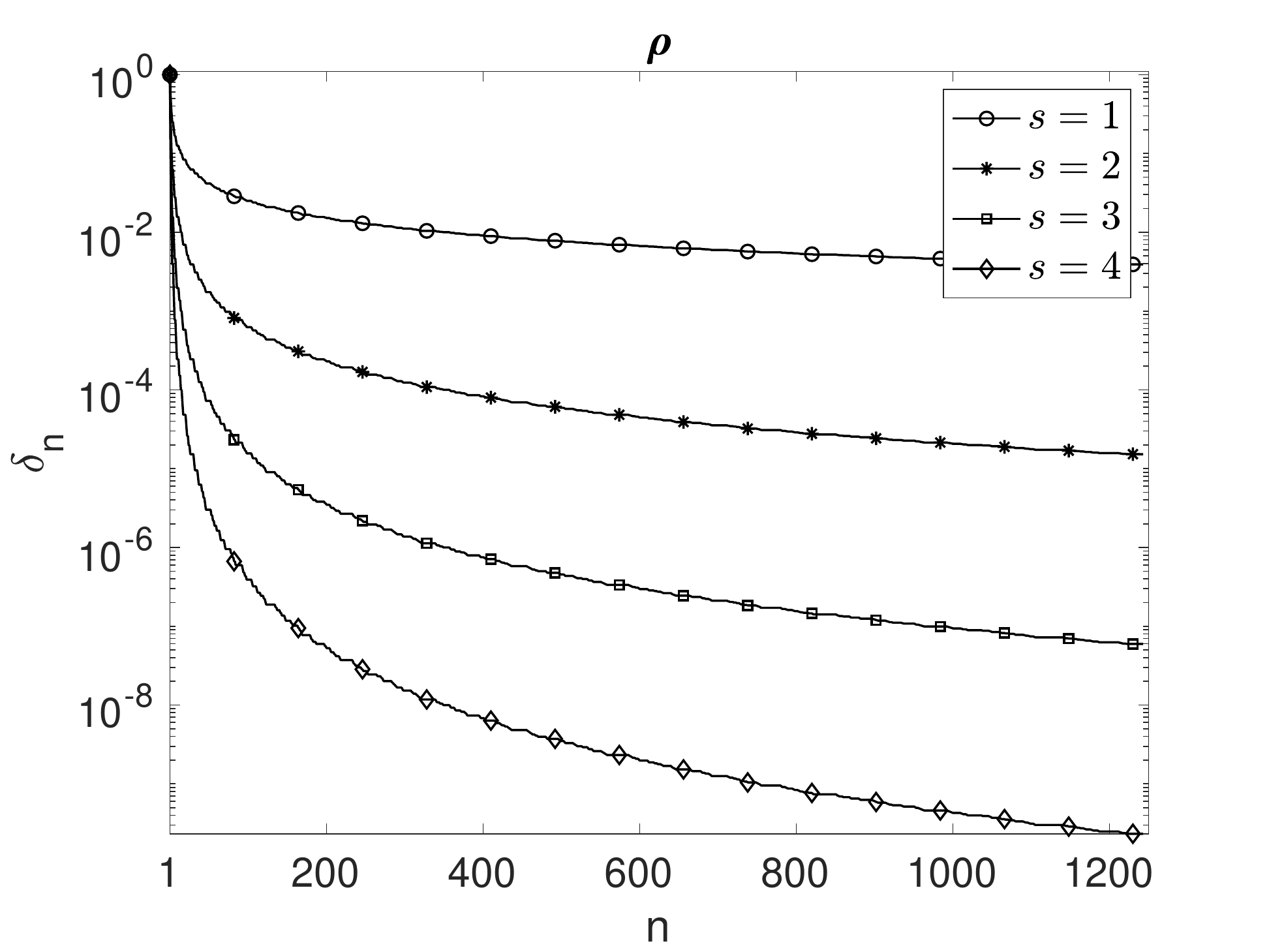}
	\caption{The graphs of $\delta_n(\rho^*(s))$ and $\delta_n(\rho(s))$ for $s=1,2,3,4$.}
	\label{fig:delta_n}
\end{figure}

In view of Remark \ref{rmk}, the behavior of the sequence $(\delta_n(\rho))_{n\ge 1}$ dictates the error of approximation for $\cB_{\rho,p}$. The values of $\delta_n(\rho(s))$ are provided in Figure \ref{fig:delta_n} for $s=1,2,3,4$ for the cases $\rho(s)=\rho^*(s)$ and $\rho(s)=((j+1)^s)_{j\ge1}$.

\subsection{The asymptotic behavior  of $\delta_n(\rho^*(s))$}
\label{sec:asymptotic}

Theorem \ref{cardtheorem}  gives an exact expression for $\#\Lambda(\e,\rho^*(s))$ which  then can be used to determine $\delta_n(\rho^*(s))$ for any $s$ and $n$. We can also use this theorem to give bounds on the asymptotic decay of $\delta_n(\rho^*(s))$. We begin with a lemma.

\begin{lemma}
\label{aslemma}
For $m=0$, $\#\Lambda( 2^{-ms},\rho^*(s))=1$, when $m=1$, $\#\Lambda( 2^{-ms},\rho^*(s))=3$, and for every $m\ge 1$, we have the following two estimates:

\begin{enumerate}[label={\rm(\roman*)}]
\item $\#\Lambda(2^{-ms},\rho^*(s))\le  2^{m+4\sqrt{m}}$,  
\item $\#\Lambda(2^{-ms},\rho^*(s))\le Cm^{-3/4} 2^{m+c\sqrt{m}}$, where $C :=(1-2^{-1/4})^{-1}$ and $ c:=\pi\sqrt{\frac{2}{3}}( \ln 2)^{-1}<4$.
\end{enumerate}
If we superimpose these inequalities we obtain
\begin{eqnarray}
\label{m2}
\#\Lambda(2^{-ms},\rho^*(s))\le 
\begin{cases}
    2^{m+4\sqrt{m}}, \quad & 2\le m\le 5, \\
    Cm^{-3/4} 2^{m+c\sqrt{m}}, & m\ge 6.
  \end{cases}
\end{eqnarray}

\end{lemma}
 
\noindent {\bf Proof:} Note that for the sequence $\rho^*(s)$ given by \eref{rhos}, we have
$$\Lambda(2^{-ms},\rho^*(s))\equiv \Lambda( 2^{-m},\rho^*(1)),\quad s>0.$$
Therefore $\Lambda(2^{-ms},\rho^*(s))$ does not depend on $s$, and in what follows we may take $s=1$. 

For the particular cases $m=0,1$ we readily check that
$$\#\Lambda( 2^{0},\rho^*(1))=1, \quad \#\Lambda( 2^{-1},\rho^*(1))=3.$$
To show  (i) and (ii), we first prove
\be
\label{l11}
\#\Lambda(2^{-m},\rho^*(1)) =\sum_{k=0}^m \#S_k \le 1+ \sum_{k=1}^m  k^{-3/4}2^{k+c\sqrt{k}}.
\ee
For {$k\ge 1$}, we note that the binomial coefficient from~\eqref{know21} can be estimated
\be
\label{energycount}
\binom {N_k-1+\#I_k}{N_k}\le  \binom {N_k +2^{k-1}}{N_k} \le  2^{kN_k},
\ee
since
\bes
\label{product}
\frac{j+2^{k-1}}{j} \le  2^k,\quad j=1,\dots,N_k.
\ees
Therefore, for any sequence  $(N_1,\dots,N_m)$ in $\cQ_m$, we have
\bes
\label{know2}
\#\Gamma(N_1,\dots,N_m)=  \prod_{k=1}^m\binom {N_k-1+ \#I_k}{N_k} \le 2^{\sum_{k=1}^mkN_k}=2^m,
\ees
yielding the estimate  
\be
\label{Smbnd}
\#S_m \le 2^m \,q(m), \quad  q(m):=\#\cQ_m.
\ee
As noted before, $q(m)$ is the same as the number of  additive partitions of the integer $m$. The number $q(m)$ has been exactly computed for small values of $m$ and there are bounds for $q(m)$ for any $m$. The following upper
bound for $q(m)$ can be found in \cite{P}:
\bes
\label{HR}
q(m)\le m^{-3/4}2^{c\sqrt{m}}, \quad m\ge 1, \quad {\rm where} \quad c=\pi\sqrt{\frac{2}{3}}( \ln 2)^{-1}.
\ees
Hence,
\bes
\label{upperbound1}
\#S_m \le m^{-3/4}2^{m+c\sqrt{m}},\quad m\ge1,
\ees
and using  Theorem \ref{cardtheorem}, we obtain \eqref{l11}.
 
We can now use \eref{l11} to prove each of the inequalities (i) and (ii). To prove (ii), it is enough to show that
\begin{equation}\label{e:induction}
1+\sum_{k=1}^m  k^{-3/4}2^{k+c\sqrt{k}}\leq Cm^{-3/4}2^{m+c\sqrt{m}}.
\end{equation}
The above relation is valid for $m=1$ and we now proceed by induction assuming that it has been proven for $m$ and verify the case $m+1$. Using the induction hypothesis, we have
\begin{align*}
  1+ \sum_{k=1}^{m+1}  k^{-3/4}2^{k+c\sqrt{k}} & \leq  (m+1)^{-3/4} 2^{m+1+c\sqrt{m+1}} + Cm^{-3/4} 2^{m+c\sqrt{m}}\\
 &   = C (m+1)^{-3/4} 2^{m+1+c\sqrt{m+1}} \left( C^{-1} + (1+1/m)^{3/4}2^{-1+c\sqrt{m}-c\sqrt{m+1}} \right)\\
 &   \leq C (m+1)^{-3/4} 2^{m+1+c\sqrt{m+1}},
\end{align*}
where to derive the last inequality we used  $\sqrt{m} < \sqrt{m+1}$, $1/m \leq 1$, and the specific value of $C$. This completes the proof of (ii).

We prove estimate (i) for $m\geq 2$ in a similar way (the case $m=1$ clearly holds) showing by induction that 
$$1+ \sum_{k=1}^m  k^{-3/4}2^{k+4\sqrt{k}} \le 2^{m+4\sqrt{m}}.$$
The details are omitted.

To prove the superimposed estimate we note that 
$$2^{m+4\sqrt{m}}\le Cm^{-3/4} 2^{m+c\sqrt{m}}\quad \text{if and only if}\quad C^{-1}m^{3/4}\le 2^{(c-4)\sqrt{m}}.$$
On the interval $[2,\infty)$, the function on the left is increasing and the function on the right is decreasing since $c<4$, and the range of $m$ for which the inequality holds is $2\le m\le 5$. The proof is completed.
\hfill $\Box$
  
\vskip .1in 
  
\begin{figure}[htbp]
\centering
\includegraphics[width=0.5\textwidth]{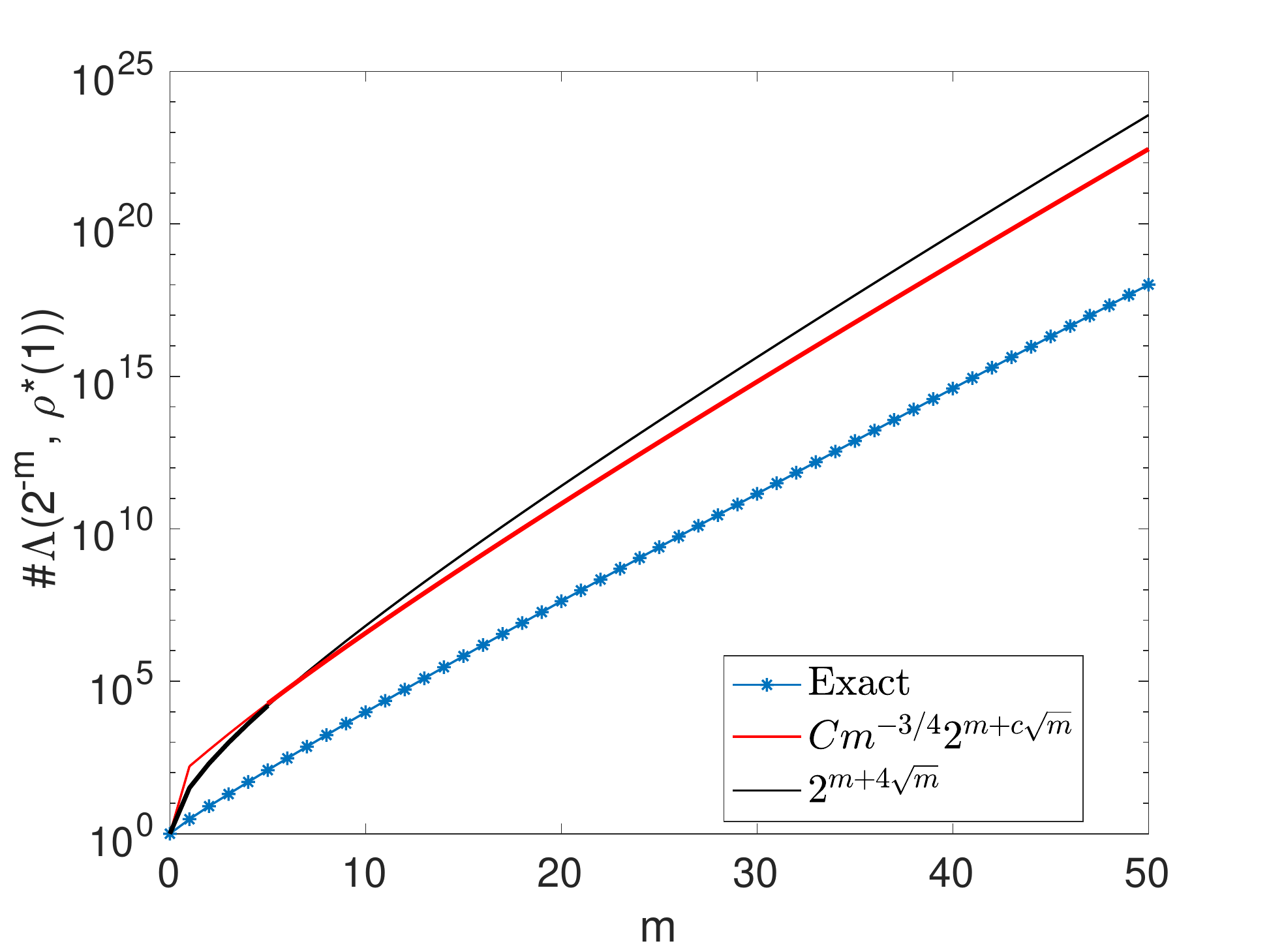}
\caption{The graphs of $\#\Lambda(2^{-m},\rho^*(1))$ and the estimates from Lemma \ref{aslemma}.}
\label{2L}
\end{figure}

In Figure \ref{2L}, we present the graphs of the exactly computed values of $\#\Lambda(2^{-m},\rho^*(1))$ compared to the estimate from Lemma \ref{aslemma}.

\subsubsection{Bounds for the error  $E_n(\cU_{\rho^*(s),p})$.}

In this section, we  use  Lemma \ref{aslemma} to give  bounds on the decay of $\delta_n(\rho^*(s))$ and $E_n(\cU_{\rho^*(s),p})$. We start with the case $p=1$.
\begin{cor} 
\label{dta0}
If $s>0$, then we have the following bounds
\be
\label{Cbound10}
\delta_n(\rho^*(s))\le 
2^{-6s}n^{-s}n^{\frac{4s\sqrt{4+\log_2 n}}{\log_2   n} }, \quad  \text{and thus}\quad 
 E_n(\cU_{\rho^*(s),1})\le 
2^{-6s}n^{-s}n^{\frac{4s\sqrt{4+\log_2 n}}{\log_2   n} },\quad n\ge 2.
\ee
\end{cor}

\noindent {\bf Proof:} We first consider the case when $n=2^k$, $k\ge 1$. Let $m$ be the largest non-negative natural number satisfying
\be
\label{satisfy1}
m+4\sqrt{m}\le k.
\ee
It follows from Lemma \ref{aslemma} that $\#\Lambda(2^{-m},\rho^*(1)) \le 2^{m+4\sqrt{m}}\le 2^k=n$. Relation \eref{dts} and the monotonicity of the sequence $(\delta_n(\rho^*(s)))_{n\geq1}$ give $\delta_n(\rho^*(s)) \le 2^{-ms}$ which, according to Remark \ref{rmk}, leads to $E_n(\cU_{\rho^*(s),1})\le 2^{-ms}$.

Let us define $\alpha$ by the equation $m=k-\alpha \sqrt{k}$ and  give an upper bound for $\alpha$. Since the integer $m+1 = k + 1 - \alpha \sqrt{k}$ does not satisfy \iref{satisfy1}, we have
\bes
k+1-\alpha \sqrt{k}+4\sqrt{k+1-\alpha\sqrt{k}}> k,
\ees
and so  
\bes
\label{need1}
16( k+1-\alpha\sqrt{k})  >   \alpha^2 k-2\alpha\sqrt{k} +1.
\ees
Rearranging terms, we have

\bes
\label{need2}
\alpha^2k+14\alpha\sqrt{k}-16k-15 < 0.
\ees
Noticing that the left-hand side vanishes for
\bes
\label{need3}
\alpha_{\pm}=-\frac{7\sqrt{k}}{k}\pm\frac{\sqrt{16k(4+k)}}{k},
\ees
we obtain the upper bound $\alpha<\alpha_+$ from which we get
\bes
\label{need4}
m=k-\alpha\sqrt{k}> k+7-4\sqrt{4+k}.
\ees
Therefore, we have the estimate 
$$2^{-ms}< 2^{-(k+7-4\sqrt{4+k})s}=2^{-7s}n^{-s}n^{\frac{4s\sqrt{4+\log_2 n}}{\log_2 n}},$$
which leads to
\bes
\label{proven}
E_n(\cU_{\rho^*(s),1})\le 2^{-7s}n^{-s}n^{\frac{4s\sqrt{4+\log_2 n}}{\log_2 n}}, \quad n=2^k.
\ees

Now, given any $n\ge 2$, we choose the largest $k$ such that $2^k\le n<2^{k+1}$. This implies that $2^{-k}<2n^{-1}$ and $\sqrt{4+k} \leq \sqrt{4+\log_2 n}$, and so we derive
\bes
\label{proved}
E_n(\cU_{\rho^*(s),1})\le E_{2^k}(\cU_{\rho^*(s),1})\le 2^{-7s}2^{-ks}2^{4s\sqrt{4+k}} \le 2^{-6s}n^{-s}2^{4s\sqrt{4+\log_2 n}}=2^{-6s}n^{-s}n^{\frac{4s\sqrt{4+\log_2 n}}{\log_2 n}},
\ees
as desired. \hfill $\Box$

The next corollary treats the case of general $p$.
\begin{cor}\label{thm:genp}
Let $1< p \leq \infty$ and let $q$ be given by $1/p + 1/q = 1$. For any $s>1/q$, we have
\be
\label{Cbound101}
E_n(\mathcal{U}_{\rho^*(s),p})\le \delta_{n,q}(\rho^*(s))\le C(q,s)n^{-s + 1/q}n^{\frac{4s\sqrt{4+\log_2 n}}{\log_2 n}},\quad n\geq 2,
\ee
where $C(q,s)$ is a constant depending only on $s$ and $q$.
\end{cor}

\noindent {\bf Proof:} The first inequality is (ii) of Theorem \ref{bestlambdatheorem}. Next, let us denote by
$$\phi(n):= n^{-s + 1/q}n^{\frac{4s\sqrt{4+\log_2 n}}{\log_2 n}}, \quad n\ge 2,$$
and observe that since $\varphi(x):=x^{\frac{4s\sqrt{4+\log_2 x}}{\log_2 x}}$ is an increasing function of $x>0$, we have
$$\phi(2^N)\le 2^{s-1/q}\phi(k), \quad  2^N< k\leq 2^{N+1}. $$
Note that  to complete the proof we need only show \eref{Cbound101} in the case $n=2^N$ because then for $2^N< k\leq 2^{N+1}$,
$$\delta_{k,q}(\rho^*(s))\leq \delta_{2^N,q}(\rho^*(s))\leq C_1\phi({2^N})\leq C_12^{s-1/q}\phi(k),\quad C_1=C_1(q,s),$$
where we have used the fact  that the sequence $(\delta_{n,q}(\rho^*(s)))_{n\ge1}$ is decreasing. Thus, we concentrate on the case $n=2^N$ and define 
\be
\label{def[si}
\psi(n):= n^{-s }n^{\frac{4s\sqrt{4+\log_2 n}}{\log_2 n}}, \quad n\ge 2.
\ee
Similarly to the function $\phi$, we have that 
$$
\psi(j)\le 2^s\psi(2^{k+1}), \quad 2^k< j\leq 2^{k+1}.
$$
It follows from Corollary \ref{dta0} that $\delta_n\le 2^{-6s}\psi(n)$, $n\ge 2$, and using the above estimate we have
\begin{eqnarray}
\label{f1}
2^{6sq}\delta_{2^N,q}^q&=&2^{6sq}\sum_{j>2^N}\delta_j^q  \le \sum_{k=N}^
\infty \sum_{2^k< j\leq 2^{k+1} } [\psi(j)]^q\le  2^{sq}\sum_{k=N}^\infty
2^k[\psi(2^{k+1})]^q\nonumber\\
&=&2^{sq-1}\sum_{k=N}^\infty2^{k+1}[\psi(2^{k+1})]^q=
2^{sq-1}\sum_{k=N}^\infty2^{-(k+1)(sq-1)}2^{4sq\sqrt{4+k+1}}\nonumber\\
&= &2^{sq-1}\sum_{k=N}^\infty2^{-k(sq-1)}2^{4sq\sqrt{4+k}}
\le 2^{sq-1}C_0^q 2^{-N(sq-1)}2^{4sq\sqrt{4+N }}.
\end{eqnarray}
Here, in the last inequality we have used the bound 
\be
\label{used}
\sum_{j=0}^
\infty2^{-j(sq-1)}2^{4sq[\sqrt{4+m+j}-\sqrt{4+m}]} \le  \sum_{j=0}^ \infty2^{-j(sq-1)}2^{4sq\sqrt{j}} \le C_0^q, \quad C_0=C_0(q,s),
\ee
valid for every $m\geq 0$, which follows from the fact that
$$\sqrt{4+m+j}-\sqrt{4+m}= \frac{j}{\sqrt{4+m+j}+\sqrt{4+m}}\le \sqrt{j}.$$
The bound \eref{f1} gives
$$
\delta_{2^N,q}\leq C_1\phi(2^{N}), \quad \hbox{with}\quad  C_1=C_1(q,s):=2^{-5s-1/q}C_0,
$$
which is \eref{Cbound101} for $n=2^N $, and therefore completes the proof of the Corollary. \hfill $\Box$

According to~\eqref{m2}, we can improve estimate \eref{Cbound10} when $n$ is large. For this, we state the following two corollaries whose proofs will be given in the appendix.
 
\begin{cor}
\label{cor:r1}
Let $m=m(n)$ be the largest natural number such that 
\be
\label{satisfy50}
\log_2 C -\frac{3}{4}\log_2 m+m+c\sqrt{m}\le \log_2 n,
\ee
where $C$ is the constant of Lemma \ref{aslemma}. Then 
\be
\label{Cbound20}
\delta_n(\rho^*(s))\le  2^{-m(n)s}, \quad \text{and therefore}\quad E_n(\cU_{\rho^*(s),1})\le 2^{-m(n)s}, \quad n\ge 2^{16}.
\ee
\end{cor}

Note that the dependence of $m$ as a function of $n$ in the above corollary is implicit. One may want to get an explicit version of that statement which is the next corollary.

\begin{cor} 
\label{cor:r2}
If  $s>0$, 
\be
 \label{Cbound2}
 \delta_{\left\lceil\tilde C\frac{n}{[\log_2 n]^{3/4}}\right\rceil}(s)\le
2^sn^{-s}
n^{\frac{cs}{\sqrt{\log_2   n}}}, \quad n\ge 2^{16},
\ee
and therefore
\be
\label{Cbound21}
 E_{\left\lceil\tilde C\frac{n}{[\log_2 n]^{3/4}}\right\rceil}(\cU_{\rho^*(s),1})\le 
2^sn^{-s}
n^{\frac{cs}{\sqrt{\log_2   n}}},\quad  n\ge 2^{16},
 \ee
where $\tilde C:=C(1-c/4)^{-3/4}$ with $C,c$ as in Lemma  \ref{aslemma}.
\end{cor}

\section{Finding the set $\Lambda(\e,\rho)$} 
\label{sec:findlambda}
 
In this section, we describe a possible strategy to build the set $\Lambda(\e,\rho)$ for any given sequence $\rho$ and a given target accuracy $\e$. A second procedure (not given here) can then be used to find $\Lambda_{n,\rho}$ when we prescribe the cardinality $n$ of the set rather than the accuracy. Before we begin describing our algorithm, let us note that   other procedures have been given for constructing $\Lambda(\e,\rho)$ (see e.g. \cite{BAS,Zech}).
 	
As above, we consider $\rho=(\rho_j)_{j\geq 1}$ to be a non-decreasing sequence such that $\rho_1>1$ and $\lim_{j\rightarrow \infty}\rho_j=\infty$. Let us denote by ${\rm supp}(\nu)$ the support of a multiindex $\nu=(\nu_1,\nu_2,\ldots)$, that is
$${\rm supp}(\nu):=\{j:\nu_j\neq 0\}.$$
Recalling the definition of $\Lambda(\e, \rho)$ given in (\ref{defLambda}), we first notice that: 
\begin{itemize}
\item $\nu=0\in \Lambda(\e, \rho)$ whenever $\e\leq 1$;
\item for every fixed $\e$, there is an index $D(\e)$ such that if $\nu\in\Lambda(\e, \rho)$, then ${\rm supp}(\nu)\subset \{1,2,\ldots,D(\e)\}$;
\item if $\e_1\leq \e_2$, then $D(\e_2)\leq D(\e_1)$.
\end{itemize}
 	
The lower set $\Lambda(\e, \rho)$ can be built using the iterative strategy described in the following Algorithm.

\begin{algorithm}
 		\caption{Construction of the lower set $\Lambda(\e, \rho)$, { $\e\leq 1$.}}\label{algo:built_lambda}
 		\begin{algorithmic}[1]
 			\BState \emph{Initialization}:
 			\State Set $T_0:=\{0\}$ and $\Lambda:=T_0$.
 			\BState \emph{Recursive Construction}:
 			\For{$i=0,1,2,\ldots$}
 			\State Set $T_{i+1}:=\emptyset$.
 			\For{$\nu\in T_i$ and $j=1,\ldots,D(\e)$}
 			\State Construct $\mu$ such that $\mu_k=\nu_k$, $k\neq j$, and $\mu_j=\nu_j+1$.
 			\If{$\rho^{\mu}\leq\e^{-1}$}
 			\State $T_{i+1} \leftarrow T_{i+1}\cup\{\mu\}$.
 			\EndIf
 			\EndFor
 			\If{$T_{i+1}=\emptyset$}
 			\State Break.
 			\Else
 			\State $\Lambda\leftarrow \Lambda\cup T_{i+1}$.
 			\EndIf
 			\EndFor
 			\Return $\Lambda$.
 		\end{algorithmic}
\end{algorithm}
 	
When implementing this algorithm in practice, we form a tree where each $\nu \in T_i$ has $D(\e)$ possible children $\mu$ to be checked for admissibility. When a constructed $\mu$ is found to be inadmissible, then it is not included in $T_{i+1}$. This stops the search down the entire subtree rooted at $\mu$. If $\mu$ is found to be admissible then it is added to $T_{i+1}$. In this way, each $T_i$ forms a level in the tree rooted with the zero sequence. When all elements of $T_i$ are exhausted, then the computation moves to processing elements in $T_{i+1}$. If the current set being processed is empty, then the procedure is ended and $\Lambda(\e, \rho) = \bigcup_{k=0}^i T_k$. Finally, we mention that the set $T_{i+1}$ corresponds to the so-called \emph{reduced margin} (see e.g. \cite{CCDS}) of the set $\bigcup_{k=0}^i T_k$.\\
 	
\begin{remark} \label{rem:optimal}
One can  deduce that the number of computations needed to construct the set $\Lambda(\e,\rho)$ is of order $\cO(m \log m)$, where $m=\#\Lambda(\e,\rho)$, provided one imposes additional growth conditions on the sequence $(\rho_j)_{j\geq 1}$ (for an analysis for another sorting algorithm see \cite{BAS}). This would cover the sequences $\rho(s)$ and $\rho^*(s)$ for example.
\end{remark}

\section{Concluding Remarks} 
\label{sec:conclusion}
 
In this work, we discussed the approximation of Banach space valued functions with an infinite number of variables by polynomials on lower sets. We defined a family of model classes $\cB_{\rho,p}$ based on anisotropic analyticity, and derived bounds for the decay rate for the approximation of these model classes using multivariate polynomials. We considered only the case when the approximation error is measured in the $L_\infty(Y,X)$ norm, though it would be interesting to develop corresponding results when measuring the approximation error in $L_q(Y,X)$ norms.  Already, several results in the case $q=2$ have been given in \cite{CM}.
 
Another setting that arises in parametric PDEs is analytic functions which have Legendre expansions (instead of Taylor expansions) with bounds on the size of the Legendre coefficients (see \cite{CCS}). It would be interesting to formally introduce and study the spaces (analogous to the $\cB_{\rho,p}$) associated to such expansions. The    functions in these spaces would now be analytic on polyellipses.
 
Our main vehicle for deriving error estimates for these classes was to use a surrogate norm in place of the $L_\infty(Y,X)$ norm. We showed in Theorem \ref{bestlambdatheorem} that for this surrogate norm, our estimates are optimal.  It would be very interesting to understand what optimal results would look like in the original $L_\infty(Y,X)$ norm, i.e., to prove lower bounds for the approximation rate in the $L_\infty(Y,X)$ norm rather than the surrogate norm.
 
We concentrated on the sequences $\rho(s)$ and $\rho^*(s)$, $s>0$, since they comply with typical assumptions in applied settings. It is possible to extend these results to more general sequences $\rho$ which eventually behave
asymptotically like $\rho(s)$ or $\rho^*(s)$. However, the behavior of the sequence in the preasymptotic regime strongly effects the final decay rate bounds for $\delta_n(\rho)$. For instance, the value $\rho_j$, representing the smoothness of $u$ in the direction $j$, might remain close to $1$ for arbitrarily many $j$ before eventually growing to $\infty$. It would be interesting to give bounds for other sequences $\rho$ with polynomial or even exponential growth.
 
Our formulation of the model classes and our approximation results have been strongly influenced by the works \cite{BCM,CM,Tran1,Zech}. The paper \cite{GO} has a significant intersection with our paper where results analogous to Corollary \ref{thm:genp} in the case $p=\infty$ are proven.

We next ellaborate on the distinctions between our paper and the results given in \cite{Tran1}. In \cite{Tran1}, the authors derive bounds for the approximation of parametric PDEs using Taylor and Legendre series. They work under the assumption that $d < \infty$, and use analyticity of the parameter-to-PDE-solution map to derive certain upper bounds on the norms of the coefficients in the Legendre and Taylor series expansions of the solution $u$. In the case of Taylor series, their analysis includes the case when $\| t_\nu \|_X \le M \rho^{-\nu}$, which corresponds to our model classes $\cB_{\rho,\infty}$. We restrict our further comments to this case. Although their results are only stated for solutions to parametric PDEs, their proofs give the following estimates for functions in $\cB_{\rho,\infty}$.
 
\begin{theorem}
\label{Trantheorem}
 
Let $\rho=(\rho_1,\dots,\rho_d) $ be a nondecreasing sequence with $\rho_1>1$. Then for any $\sigma > 0$, there exists an $n(d,\sigma)$ such that for all $n\ge n(d,\sigma)$, 
\be
\label{Tranbound}
E_n(\cU_{\rho,\infty})\le  C_\sigma n \exp\left(-\left( \prod_{i=1}^d \log(\rho_i) \frac{n d!}{1+\sigma} \right)^{1/d} \right),
\ee
holds with $C_\sigma:= (4e+4\sigma e -2)\frac{e}{e-1}$.
\end{theorem}
If we specialize to the sequence $\rho^*(s)$, $s>0$,  then their result takes the form
\bes
\label{TWrate}
E_n(\cU_{\rho,\infty})\le Cne^{-c(d,s)n^{1/d}}, \quad n\ge n(d,\sigma),
\ees
where $C$ has an absolute bound and $c(d,s)$ actually grows with $d$ and $s$. Note that the bound is subexponential in $n$, and hence is better than the algebraic rate given in our estimates. The reason for this is the assumption that $d$ is finite. However, we must emphasize that  the number $n(d,\sigma)$ grows exponentially in $d$, and so this result can only be applied when $n$ is very large. We have concentrated on obtaining results that hold for all $n$ and all $d$ with no dependence on $d$.
	
The reason for this restriction on $n$  in \cite{Tran1} is  that their proof of this theorem utilizes bounds on the number of lattice points  $t\N^d$ in the simplex $S=S(1/\ln \rho_1,\dots,1/\ln \rho_d)$. Their bound requires that this number behaves like $t^{-d}\meas(S)$. As discussed in the remarks following the proof of Lemma \ref{countlemma}, this asymptotic count on the lattice points is effective only for $t$ small and in turn $n$ prohibitively large.
 
By contrast, our results given above apply for $d=\infty$ and any $n$. When $d$ is finite we can always extend the sequence to an infinite sequence in an arbitrary way. In this way our results apply without any restrictions on the size of $n$ relative to $d$.

\section{Appendix:  Proofs of Corollaries \ref{cor:r1} and \ref{cor:r2}}
\label{appendix}

\noindent {\bf Proof of Corollary \ref{cor:r1}}:
Let $m=m(n)$ be the largest natural number satisfying (\ref{satisfy50}). One can check that for $n\ge 2^{16}$, we have $m(n)\ge 6$, and thus it follows from (ii) or Lemma \ref{aslemma} that
$$\#\Lambda(2^{-m(n)},\rho^*(1)) \le Cm(n)^{-3/4}2^{m(n)+c\sqrt{m(n)}}\le n,$$
which gives $\delta_n(\rho^*(s)) \le 2^{-m(n)s}$, and thus $E_n(\cU_{\rho^*(s),1}) \le 2^{-m(n)s}$. 
\hfill $\Box$

\vskip .1in
\noindent {\bf Proof of Corollary \ref{cor:r2}}: To show \eref{Cbound2}, we proceed as follows. We consider first the case $n=2^k$, $k\ge 16$. Let $m$ be the largest non-negative natural number satisfying
$$m+c\sqrt{m}\le k,$$
and let $\beta$ be defined by the equation $m=k-\beta\sqrt{k}$. Since $k\ge 16$, the largest $m$ that satisfies the above estimate is greater or equal to $6$. Moreover, we can easily show that $\beta\le c$. Therefore, we use the fact that $m=k-\beta\sqrt{k}\ge k-c\sqrt{k}$ and that $k-c\sqrt{k}\ge (1-c/4) k$ for $k\ge 16$, which gives
$$\log_2 m\ge \log_2 (1-c/4) +\log_2 k.$$
Thus if $C_1:=\log_2 C$, we have
\be
\label
{satisfy41}
C_1 -\frac{3}{4}\log_2 m+m+c\sqrt{m}\le 
C_2-\frac{3}{4}\log_2k +k, \quad C_2:=C_1 -\frac{3}{4}\log_2 (1-c/4).
\ee
It follows (since $m\ge 6$) that
$$\#\Lambda(2^{-m},\rho^*(1))\le Cm^{-3/4}2^{m+c\sqrt{m}}\le\tilde Ck^{-3/4}2^k
=\tilde C\frac{n}{[\log_2 n]^{3/4}}, \quad \tilde C:=C(1-c/4)^{-3/4}>1.$$
Therefore, \eref{dts} and the monotonicity of the sequence $(\delta_n(\rho^*(s)))_{n\ge 1}$ give 
$$\delta_{\left\lceil\tilde C\frac{n}{[\log_2 n]^{3/4}}\right\rceil}(\rho^*(s))\le2^{-ms}\le n^{-s}n^{\frac{c s}{\sqrt{\log_2   n}} },\quad n= 2^{k}, \quad k\ge 16,$$ 
which, according to  Remark \ref{rmk} leads to
$$E_{\left\lceil\tilde C\frac{n}{[\log_2 n]^{3/4}}\right\rceil}(\cU_{\rho^*(s),1})\le n^{-s}n^{\frac{c s}{\sqrt{\log_2   n}} }.$$
Now, if $k\ge 16$ is such that $2^k\le n<2^{k+1}$, it follows that 
\bes
E_{\left\lceil\tilde C\frac{n}{[\log_2 n]^{3/4}}\right\rceil}(\cU_{\rho^*(s),1})\le 
E_{\left\lceil\tilde C\frac{2^k}{k^{3/4}}\right\rceil}(\cU_{\rho^*(s),1})\le 2^{-ks } 2^{cs \sqrt{k} }\le 2^sn^{-s} (2^{\log_2 n})^{\frac{c{ s}}{\sqrt{\log_2 n}}}
{ =} 2^sn^{-s}
n^{\frac{cs}{\sqrt{\log_2   n}}},
\ees
which is \eref{Cbound2}.
\hfill $\Box$

\vskip .2in
\noindent {\bf Acknowledgements} The authors would like to acknowledge and thank Matthew Hielsberg for the help in carrying the numerical experiments.

\bibliographystyle{plain}

\vskip .25in

\noindent
\{AB, RD, DG, PJ, GP\},  Dept. of Mathematics. Texas A$\&$M University.  College Station, TX, 77843,  \{bonito,rdevore,dguignard,pjantsch,gpetrova\}@math.tamu.edu

\end{document}